\documentclass[10pt,leqno]{amsart}
\usepackage{graphicx}
\usepackage{indentfirst,csquotes}

\usepackage{amssymb,amsthm,amsmath}
\usepackage{xcolor,paralist,fancyhdr,etoolbox}
\usepackage{hyperref}
\newtheorem{theorem}{Theorem}[]
\newtheorem{definition}[theorem]{Definition}

\newtheorem{lemma}[theorem]{Lemma}
\newtheorem{proposition}[theorem]{Proposition}
\newtheorem{corollary}[theorem]{Corollary}

\newtheorem{remark}[theorem]{Remark}
\newtheorem{assumptions}{Assumptions}

\hypersetup{ colorlinks=true, linkcolor=blue, filecolor=black, urlcolor=black }

\newcommand{\Prob}[2][]{\mathbb{P}{#1} \left[ #2 \right]}

\newcommand{\Esp}[2][]{\mathbb{E}{#1} \left[ #2 \right]}

\newcommand{\Tr}[1]{\operatorname{Tr} \left( #1 \right)}

\newcommand{\diag}[1]{\operatorname{diag} \left( #1 \right)}

\newcommand{\Abs}[1]{\left| #1 \right|}

\newcommand{\Norm}[1]{\left\lVert #1 \right\rVert}

\newcommand{\Braces}[1]{\left\lbrace #1 \right\rbrace}

\newcommand{\Parentheses}[1]{\left( #1 \right)}

\newcommand{\Angles}[1]{\left\langle #1 \right\rangle}

\newcommand{\Ind}[1]{\mathbb{1}_{ #1 } }

\newcommand{\diff}[0]{\operatorname{d}}

\newcommand{\cadlag}{\emph{càdlàg}}

\newcommand{\EE}[0]{\mathbb{E}}
\newcommand{\FF}[0]{\mathbb{F}}

\newcommand{\HH}[0]{\mathbb{H}}

\newcommand{\NN}[0]{\mathbb{N}}
\newcommand{\PP}[0]{\mathbb{P}}

\newcommand{\RR}[0]{\mathbb{R}}

\newcommand{\Bb}[0]{\mathcal{B}}

\newcommand{\Dd}[0]{\mathcal{D}}
\newcommand{\Ee}[0]{\mathcal{E}}
\newcommand{\Ff}[0]{\mathcal{F}}
\newcommand{\Gg}[0]{\mathcal{G}}
\newcommand{\Hh}[0]{\mathcal{H}}

\newcommand{\Ll}[0]{\mathcal{L}}
\newcommand{\Mm}[0]{\mathcal{M}}
\newcommand{\Nn}[0]{\mathcal{N}}

\newcommand{\Pp}[0]{\mathcal{P}}

\newcommand{\Rr}[0]{\mathcal{R}}
\newcommand{\Ss}[0]{\mathcal{S}}

\newcommand{\Xx}[0]{\mathcal{X}}

\begin{document}


\title[Conditional McKean-Vlasov SDEs]{Conditional McKean-Vlasov Differential Equations with Common Poissonian Noise: Propagation of Chaos} 
\date{\today}
\thanks{This work was supported by the National Council of Science and Technology (CONACyT), scholarship number 863210, as part of the PhD Thesis of the second author.}


\author{Daniel~Hernández-Hernández}
\address{Centro de Investigación en Matemáticas A.C. Calle Jalisco s/n. 36240 Guanajuato, México}
\email{dher@cimat.mx}

\author{Joshué Helí~Ricalde-Guerrero}
\address{Centro de Investigación en Matemáticas A.C. Calle Jalisco s/n. 36240 Guanajuato, México}
\email{joshue.ricalde@cimat.mx}


\begin{abstract}
A model for the evolution of a large population  interacting system is considered in which a marked Poisson processes influences their evolution, together with a Brownian motion. Mean field McKean-Vlasov limits of such system are formulated  studying first both systems individually. Letting the population size growing to infinite, the weak convergence of the solutions of such systems is proved; in other words, propagation of chaos of such systems is obtained.

\smallskip
\noindent \textsc{Keywords:} Propagation of chaos; Interacting particle systems; Mean field limit; Random environment; McKean-Vlasov system.

\smallskip
\noindent \textsc{MSC2020 Classification:} 35Q70; 60B10; 60K35; 65C35; 82C22.

\end{abstract} 


\maketitle



\section{Introduction}
\label{Sec:Introdution}
In this paper we are concerned with the following model: For a probability space $(\Omega,\Ff,\PP)$  equipped with a general filtration $\FF=\{\Ff_t\}$, consider the class of $\RR^d$-valued processes $X^{n} = \Braces{ X^{i,n}, i=1,\ldots,n }$, defined as the solution to the \textit{finite-population} stochastic differential equation (SDE)
\begin{align}
    \label{Eq:X^in}
    X^{i,n}_t
    =&
    X^{i,n}_0
    +
    \int_0^t b_s( \overline{ \mu }^{n}_{s-}, X^{i,n}_{s-} ) \diff s
    +
    \int_0^t
        \sigma_s( \overline{ \mu }^{n}_{s-}, X^{i,n}_{s-} )
    \diff W_s^{i}
    \\ \nonumber
    &+
    \int_{[0,t] \times \mathbf{R}} \gamma_s( \overline{ \mu }^{n}_{s-}, r, X^{i,n}_{s-})
    \tilde{\eta}^{\lambda,n}(\diff s, \diff r),\;\;\;t\leq T,
\end{align}
for  $i=1,\ldots,n$, where $\overline{ \mu }^{n}$  denotes the empirical measure process 
\begin{align}
    \label{Eq:Overline_mu^n}
    \overline{ \mu }^{n}
    &=
    \Braces{ \overline{ \mu }^{n}_t = \frac{1}{n} \sum_{i=1}^n \delta_{ X^{i,n}_t },\ t \in [0,T] }.
\end{align} 
Given an arbitrary  finite collection of indexes $I \subset \{ 1, \ldots, n \}$, the present paper analyzes the (weak) convergence of the subclass $X^{I,n} = \Braces{ X^{i,n}, i \in I } \subset X^n$ to a collection of pathwise-exchangeable processes $X^I = \Braces{ X^{i}, i \in I }$ when $n \to \infty$, such that each $X^i$ solves the \emph{conditional McKean-Vlasov} SDE with marked jumps
\begin{align}
    \label{Eq:Underline_X}
    X^{i}_t
    =&
    X^{i}_0
    +
    \int_0^t b_s( \Ll^1\Parentheses{ X^{i} }_{s-}, X^{i}_{s-} ) \diff s
    +
    \int_0^t
        \sigma_s(\Ll^1\Parentheses{ X^{i} }_{s-}, X^{i}_{s-} )
    \diff W_s^{i}
    \\ \nonumber
    &+
    \int_{[0,t] \times \mathbf{R}} 
        \gamma_s( \Ll^1\Parentheses{ X^{i} }_{s-}, r, X^{i}_{s-} )
    \tilde{\eta}^{\lambda}(\diff s, \diff r),
\end{align}
with $\Ll^1( X^i )$  representing the conditional measure process
\begin{align}
    \label{Eq:Ll}
    \Ll^1( X^i )
    &:=
    \Braces{ 
        \Ll^1( X^i )_t 
        := 
        \Ll\Parentheses{ X^i \middle\vert \Ff_T^{ (\lambda, N) } }_t, 
    \ t \in [0,T] }.
\end{align}
In \eqref{Eq:X^in} and \eqref{Eq:Underline_X}, $X^{i,n}_0$ and $X^i_0$ are independent, identically distributed copies of a random variable $X_0 \in L^2(\RR^d )$ for all $n\in \NN$ and each $i\in \{1,\ldots,n\}$, and $W = \Braces{ W^i, i=0,1\ldots }$ is a sequence of $k$-dimensional, independent Brownian motions; $\eta^{\lambda,n}$ is the lifting of an $l$-dimensional marked Poisson process $(N^{\lambda,n},\xi^{\lambda,n})$ on $[0,T]$ and marks on a Polish space $(\mathbf{R},\Rr)$, with an $\FF$-admissible kernel $K^{\lambda,n}_t(\diff r)$ and compensated martingale measure
\begin{align*}
    \tilde{\eta}^{\lambda,n}(\diff t, \diff r) 
    :=
    \eta^{\lambda,n}(\diff t, \diff r) - K^{\lambda,n}_t(\diff r) \diff t.
\end{align*}
For $\tilde{\eta}^{\lambda}$ in (\ref{Eq:Underline_X}) there is a similar description, in terms of  $(N^{\lambda},\xi)$ and $K^{\lambda}_t(\diff r)$ \cite{bremaud_point_2020}.

Throughout  we assume without loss of generality that $\eta^{\lambda,n}$ and $\eta^\lambda$ are characterised by an independent Poisson process $N$ on $[0,T] \times \mathbf{R} \times \RR_+$ with intensity measure $\diff t \otimes Q(\diff r) \otimes \diff s$ (to be specified later) for some probability measure $Q$ on $(\mathbf{R},\Rr)$. For a given  non-negative, predictable, locally integrable process $\lambda$, we assume that  kernels $K^{\lambda,n}$ and $K^{\lambda}$ can be decomposed as
\begin{align*}
     &K_t^{\lambda,n}(\diff r) = \lambda_t(\overline{\mu}^n_{t-}, r ) Q(\diff r),
     &
     &\mbox{and}
     &
     &K^{\lambda}_t(\diff r) = \lambda_t( \Ll^1(X)_{t-}, r ) Q(\diff r),
\end{align*}

Inspired by the microscopic models in thermodynamics and statistical mechanics developed in the last century \cite{brush_further_2003}, the theory of \emph{propagation of chaos} was born from a rigorous formulation for the study of large systems of interacting particles by M. Kac, in his 1956's seminal work \cite{kac_foundations_1956}, where he mentions his  ``motivation by the mathematical justification of the classical collisional kinetic theory of Boltzmann for which he developed a simplified probabilistic model'' \cite{chaintron_propagation_2022}. 

To understand the notion of \emph{chaos} introduced by Kac, first consider a population of symmetric, or \emph{exchangeable}, particles; i.e. a collection of random elements satisfying that  the joint distribution of any finite sample from the population is invariant under each permutation of the sampled particles. Then, the population is said to be \emph{chaotic} when it holds that all the particles become statistically independent from each other as the size of the population goes to infinity. As a result, a dynamical system that \emph{propagates chaos} can be heuristically interpreted as a system where every state will be chaotic if the initial conditions had such characteristic. More recently, several studies \cite{bodineau_brownian_2016}, \cite{gallagher_newton_2014}, \cite{mischler_kacs_2013} on the propagation of chaos have proposed interesting responses to questions raised decades ago in Boltzman's original theory, adapting novel models to problems in various disciplines, such as biology, medicine, computer science and social science, to be explained in the context of mathematical kinetic theory (\cite{muntean_collective_2014}, \cite{albi_vehicular_2019}, \cite{degond_mathematical_2019}, just to name a few recent examples). Moreover, the results on propagation of chaos have also played an important  role in the development of the theory of \emph{mean-field games} and \emph{mean-field optimal control}. In both cases,  optimization problems involving a very large number of interacting agents with similar goals are considered, where the individual contribution to the overall system becomes negligible as the population size increases  \cite{lasry_jeux_2006-1}, \cite{lasry_jeux_2006}, \cite{caines_large_2006}, \cite{carmona_mean_2016}. 

Historically, the study of the conditional propagation of chaos  for discontinuous noise processes has been less studied than its continuous counterpart; to our knowledge, the most general result in this direction comes from mean-field theory, where a Gaussian white noise field is taken into consideration (see \cite[Chapter 2]{carmona_probabilistic_2018-1}). While  continuous (Gaussian) shock models do provide an incredibly flexible and relatively easy-to-use set of tools for many situations (at least in comparison with other  discontinuous processes), jumps are by their  nature one of the best suited ways to model atypical and unexpected events. This fact of course has not gone unnoticed, as many recent studies -most of them within the framework of mean-field theory- have indeed considered the inclusion of jumps into their models; see, for instance,  \cite{shen_maximum_2013}, \cite{hafayed_mean-field_2014}, \cite{bensoussan_mean-field-type_2020}, \cite{benazzoli__2019}, \cite{zhang_general_2018}, \cite{ledger_at_2021}. 

Besides some recent advances on the study of diffusion processes, the main feature and  novelty of the  model proposed  in this paper is that, unlike the usual continuous common noise, the population is affected by discontinuous perturbations that appear at random time intervals.  Physically, our model describes the behavior of a large population of particles inside a box that individually follow a diffusion with certain parameters and only interact with each other through the average position of the whole population; then, at random times an alarm sounds (represented here by the jump times of the marked process), and the box is shaken in a specific way (the corresponding mark), causing  a change in the position of the particles. Moreover, the frequency with which this alarm occurs increases as the variability of the population increases, reflecting a stabilizing behavior close to an equilibrium.

The system of SDEs (\ref{Eq:X^in})-(\ref{Eq:Underline_X}) studied in this paper   is inspired by  well-known problems of \emph{systemic risk} from economy and finance \cite{fouque_handbook_2013}. The game theoretic formulation described here allow us to interpret each particle as a \emph{bank} (\emph{agent}) within a \emph{banking financial market}. Then, the average position of the particles serves as a benchmark for the inter-bank lending activities in the market.  In this context, the proposed model admits the presence of \emph{black swans}; that is, uncommon or rare events that do not occur frequently, but affect the entire market, and are likely to have long-lasting effects due to the presence of a self-excitation component. The marks in this case can be interpreted as a way of representing the strength or scale with which the event will affect the population.

The main result of this work, Theorem \ref{Thm:Propagation-of-Chaos}, establishes that as the population size goes to infinity, the weak limit of the law of any finite set of indexes corresponds to   the mean-field McKean-Vlasov SDE, also called as  conditional propagation of chaos. This result is presented in Section \ref{Sec:Propagation}, together with some preliminary results needed for its proof. The rest of the paper  is structured as follows. In Section \ref{Sec:Prelim}, we introduce the probabilistic framework under which the stochastic process evolves, including the noise processes. The existence of strong solutions for the finite population interacting system is analyzed in Section  \ref{ExistenceStrongSol}. We highlight that our methodology allows the inclusion of \emph{Regime-switching diffusions} \cite{shao_strong_2015}  \cite{shao_propagation_2022}, another prominent type of dynamical systems, which can be  explained in the context of a common Poissonian noise; the results applied to this type of systems are presented in Section \ref{Sec:(Example)Regime-switching-diffusions}.  Finally, for the sake of readability, we include in Appendix \ref{Sec:Canonical-space} a digression on the conditional law as a random environment, which is necessary to present the arguments involved in the proof of the main result.

\section{The elements of the model and preliminaries}
\label{Sec:Prelim}

Throughout this paper, we use $(x\cdot y) = x^\top y$ to denote the usual inner product in $\RR^d$; if $x \in \RR^d$, we write $x^{-i} = \Parentheses{ x^1, \ldots, x^{i-1}, x^{i+1}, \ldots, x^d }$, $i=1, \ldots, d$ to represent the vector $x$ without $i$-th coordinate. For matrices, we use the Frobenius inner product, i.e. $(A \cdot B)_\mathrm{Fr} = \Tr{A B^\top}$ for $A,B$ real-valued matrices of the same dimensions, with the corresponding norm given by $\Norm{A}_\mathrm{Fr}^2 = \Tr{ A A^\top }$.

Given the time interval $[0,T]$, with $T>0$ fixed, and a Polish space $(\Xx,d_\Xx)$, we denote by 
\begin{align*}
    \mathrm{C}([0,T];\Xx)
    :=&
    \Braces{X:[0,T] \to \Xx \ \middle\vert \ X\mbox{ has continuous trajectories }},
    \\
    \mathrm{D}([0,T];\Xx)
    :=&
    \Braces{X:[0,T] \to \Xx \ \middle\vert \ X\mbox{ has \cadlag\ trajectories }}.
\end{align*}
In particular, we write $\Dd([0,T];\Xx)=\Bb\Parentheses{ \mathrm{D}([0,T];\Xx) }$ for the Borel $\sigma$-algebra  with respect to the Skorohod topology \cite{billingsley_convergence_2013}. We also denote the space of measures $\Mm$ and $\Mm_c^*$ as
\begin{align*}
    \Mm(\Xx) :=& \Braces{ \nu \middle\vert \nu\mbox{ is a  finite measure over }(\Xx,\Bb(\Xx)) },
    \\
    \Mm_c^*(\Xx) :=& 
        \Braces{ \nu \in \Mm(\Xx) \middle\vert \nu\mbox{ is a simple, counting measure over }(\Xx,\Bb(\Xx)) }.
\end{align*}

Denote by $\Pp(\Xx)$ the set of Borel probability measures on $\Xx$ and, for $p \in [1,\infty)$, define the set of probability measures with finite moments of order $p$ as
\begin{align*}
    \Pp_p(\Xx)
    := 
    \Bigg\lbrace 
        \nu \in \Pp(\Xx) 
        \ :\ 
        \int_\Xx d_\Xx(x_0,x)^p \nu( \diff x ) < +\infty \  \mbox{for  any } x_0 \in \Xx
    \Bigg\rbrace.
\end{align*}
Given $\nu^1,\; \nu^2\in \Pp(\Xx) $, define $\Pi(\nu^1,\nu^2) \subset \Pp(\Xx \times \Xx)$ as the set of probability measures over the product space $\Xx \times \Xx$ with marginals $\nu^1$ and $\nu^2$, respectively. The $p$-th Wasserstein distance between two probability measures $\nu^1, \nu^2 \in \Pp_p(\Xx)$ is then defined as
\begin{equation}
    \label{Eq:W_p^p}
    W_p (\nu^1, \nu^2)^p
    :=
    \min \Braces{ 
        \int_{ \Xx \times \Xx } d_\Xx( x_1, x_2 )^p \pi( \diff x_1, \diff x_2 )
        \ :\ 
        \pi \in \Pi( \nu^1, \nu^2 )
    }.
\end{equation}
Throughout we will be making use of the following inequality when $p=2$,
\begin{align}
    \label{Eq:W_2-Inequality}
    W_2 (\nu^1, \nu^2)^2 \leq \Esp{ d_\Xx(Y^1,Y^2)^2 },
\end{align}
where $Y^1,Y^2$ are independent $\Xx$-valued random variables defined on the same probability space $(\Omega,\Ff,\PP)$ and $\nu^i = \PP \circ (Y^{i})^{-1}$ for $i=1,2$.

Given $\Theta$ a stochastic process, we will some times need to introduce other filtrations associated with $\Theta$. The natural filtration is denoted by $\FF^{ \mathrm{nat},\Theta }$, while the $\PP$-null complete right-continuous augmentation is denoted by $\FF^{ \Theta }$. In general, we say that a filtration $\FF$ satisfies the usual conditions if it is right-continuous and complete.

Any mapping $\Psi:\Omega \times [0,T] \times \Pp_2(\RR^d) \times \RR^d \to \RR^p$ is denoted by $\Psi_t(\nu,x)(\omega)$ for $(\nu,x) \in \Pp_2(\RR^d) \times \RR^d$ and $(\omega, t) \in \Omega \times [0,T]$. Furthermore, for each $(\nu,x),(\nu',x')\in \Pp_2(\RR^d) \times \RR^d$, its difference is denoted by $ \delta \Psi_t(\nu,x)(\omega)= \Psi\Parentheses{ \omega, t, \nu, x } - \Psi\Parentheses{ \omega, t, \nu', x' }$. This notation is also used to represent the difference of any two elements of an Euclidean space $x,x'$  as $  \delta x = x - x'$.

{\bf Spaces for stochastic processes.} 
We write $\Hh^2(\RR^d)$ (or simply $\Hh^2$ when there is no ambiguity) for the Hilbert space of $\RR^d$-valued, square-integrable martingales (with respect to $\FF$), i.e.
\begin{align*}
    \Hh^2(\RR^d)
    :=
    \Bigg\lbrace 
		M:\Omega \times &[0,T] \to \RR^d
		\ \Bigg\vert\ 
        M\mbox{ is an }\FF\mbox{-martingale with } \sup_{ [0,T] }\Esp{ \Abs{ X_t }^2 }<\infty
    \Bigg\rbrace.
\end{align*}
This space is endowed with  the norm $\Norm{\ \cdot\ }_{\Hh^2(\RR^d)}$, induced by  the inner product $\Parentheses{ M \cdot M' }_{\Hh^2}=\Esp{ M_T M'_T }$, defined  by $\Norm{ M }_{ \Hh^2 }:=\Esp{ \Abs{ M_T }^2 }^{ \frac{1}{2} }$.  For any subset $B \subset \Hh^2$, $B^\perp \subset \Hh^2$ denotes the orthogonal set of $B$, i.e. the set of elements $M\in \Hh^2$ being  orthogonal to some martingale in $B$. The space of well-posed integrands with respect to $M$  is denoted by $\HH^2(M)$, and it is endowed with the norm of predictable quadratic variation $\Norm{\ \cdot\ }_{\HH^2(M)}$. If $M'$ is another square integrable martingale, we write $\HH^2(M,M') = \HH^2(M) \cap \HH^2(M')$. We denote by $\Ll^2(M)$ the space of martingale-driven, Itô stochastic integrals of processes in $\HH^2(M)$ with respect to $M$, i.e.  $\Ll^2(M):=\Braces{ \int H \diff M\ :\ H \in \HH^2(M) },$ equipped with the martingale norm \cite[Ch. IV.2]{protter_stochastic_2005}.  We use the notation $[H,H']$ to denote the \textit{quadratic covariation} between semimartingales $H$ and $H'$, and $\Angles{H,H'}$ to denote its compensator when it exists; this process is also known as \textit{conditional} or \textit{predictable quadratic covariation} \cite[Ch. III.5]{protter_stochastic_2005}. In particular,  we write $[H]=[H,H]$ and $\Angles{H} = \Angles{H,H}$. The space of adapted, $\Xx$-valued \cadlag\ processes with finite supremum norm is denoted by $(\Ss^2(\Xx), \Norm{\ \cdot\ }_\ast$). In the $d$-dimensional case,
\begin{align*}
    &\Norm{ M }_*^2 = \max_{1 \leq i \leq d} \Norm{ M^{(i)} }_*^2 
        = \max_{1 \leq i \leq n} \Esp{ \sup_{[0,T]} \Abs{ M^{(i)}_t }^2 },
\end{align*}
for all $M = (M^{(1)},\ldots,M^{(d)}) \in \Ss^2(\Xx)$. Moreover, $\Parentheses{ \Ss^2(\Xx), \Norm{\ \cdot\ }_* }$ is a Banach space \cite[Ch. V.2]{protter_stochastic_2005}. We use the symbol $\FF' \hookrightarrow \FF$ whenever $\FF'$ and $\FF$ are filtrations on the $\sigma$-algebra $\Ff$ such that $\FF' \subset \FF$ and  $\FF'$-martingales are also $\FF$-martingales. In particular, an $\FF$-adapted \cadlag\ process $\Theta = \Parentheses{ \Theta_t, 0 \leq t \leq T }$ with values in a Polish space is said to be \emph{compatible with $\FF$ (under $\PP$)} if $\FF^\Theta\hookrightarrow\FF$. Product spaces will always be endowed with the product $\sigma-$field.

\noindent
{\bf Random environment.}  Let $\mu$ be a square-integrable process $\mu$ taking values in the Wasserstein space $(\Pp_2(\RR^d),W_2)$. In the system under study, this process plays the role of an endogenous random environment influencing the dynamics of the SDEs. More precisely, let $\mu$ be a measurable map from $(\Omega,\Ff,\PP)$ into 
\begin{align*}
    \Parentheses{ 
		\mathrm{D} \Parentheses{ [0,T]; \Pp_2(\RR^d) }, 
		\Dd \Parentheses{ [0,T]; \Pp_2(\RR^d) } 
    }.
\end{align*}
For each $0 \leq t \leq T$, denote by $\mu_t$  the projection of $\mu$ at time $t$. We say that $\mu$ is a \emph{random environment }if there exists an $\Ff_0^{\mu}$-measurable variable $X_0$ and a probability measure $\nu_0 \in \Pp_2(\RR^d)$,  the \emph{initial distribution of $\mu$}, such that
\begin{align}
	\label{Eq:Conditions-on-Random-Environment-I}
	\Ll(X_0) &= \nu_0,
	&
	\Prob{  \mu_0 = \nu_0} &= 1,
\end{align}
and, for each $\nu\in\Pp_2(\RR^d)$,
\begin{align}
	\label{Eq:Conditions-on-Random-Environment-II}
	\Esp{ \sup_{t\in [0,T]} W_2(\mu_t,\nu)^2 } < \infty.
\end{align}
In our context, $\mu$ concentrates along the time the states of the entire set of interacting particles into a single population statistic.

\textbf{Admissible set-up.} In this work, we are interested in \emph{marked Poisson processes} with random intensity and marks on a measurable space $(\mathbf{R},\Rr)$, such that the information from the environment influences their behavior. This is done by assuming first that $\FF^{\mu} \subset \FF$,  and that $\lambda : \Omega \times [0,T] \times \mathbf{R} \to \RR^l_+$ is an \textit{$\FF$-admissible intensity candidate}; that will mean that $\lambda$ satisfies the following conditions.
\begin{enumerate}
    
    \item
        $\lambda_t^{(i)}(\cdot,r) \perp \lambda_t^{(j)}(\cdot,r)$ for $i \neq j$, and each $(t,r) \in [0,T] \times \mathbf{R}$;
    
    \item 
        $\lambda$ is $\FF$-predictable;
    
    \item
        for all $0 \leq a < b \leq T$,
        \begin{align*}
            &\int_a^b \int_{\mathbf{R}} \lambda_t^{(i)}(r) Q^{(i)}(\diff r) \diff t < \infty
            &
            &\PP-\mbox{a.s.},
    \end{align*}
    where $Q^{(1)},\ldots,Q^{(l)}$ are  finite measures on $(\mathbf{R},\Rr)$.
    
\end{enumerate}
Then, we say that the random kernel $K^\lambda : \Omega \times [0,T] \times \Rr \to \RR_+^l$, defined entry-wise as
\begin{align}\label{KerK}
    &K^{\lambda(i)}_t(\omega, C) := \int_{C} \lambda^{(i)}_t(\omega,r) Q^{(i)}(\diff r),
    &
    &i = 1, \ldots, l,
\end{align}
for all $(\omega,t) \in \Omega \times [0,T]$ and all $C \in \Rr$, is an \emph{admissible $\FF$-intensity kernel}. 

\begin{remark}
In the context of point processes, $K^\lambda$ is simply an $\FF$-predictable and locally integrable kernel from $(\Omega \times [0,T],\Ff \otimes \Bb([0,T])$ into $(\mathbf{R},\Rr)$; i.e.,
\begin{itemize}
    \item 
        $K^\lambda_\cdot(\cdot,C)$ is $\FF$-predictable for each $C \in \Rr$,
    \item 
        $K^\lambda_t(\omega,\cdot)$ is a finite measure on $(\mathbf{R},\Rr)$ for each $(\omega,t) \in \Omega \times [0,T]$.
\end{itemize}
\end{remark}

The \textit{admissible noise process} is defined as a marked Poisson process $(N^\lambda,\xi)$ on $[0,T]$ with marks on $\mathbf{R}$, with $\FF$-intensity kernel $K^\lambda$. This means that $N^\lambda$ is a Poisson process with stochastic intensity $K^\lambda_t(\omega,\mathbf{R}) = \int_\mathbf{R} \lambda_t(\omega,r)Q(\diff r)$ and jumping times $\tau = \Braces{\tau_n, n \geq 1 }$:
\begin{align*}
    \tau_{n + 1}(\omega)
    = 
    \inf\Braces{ t \geq \tau_{n}(\omega) ; \int_0^{t} K^\lambda_s(\omega,\mathbf{R}) \diff s \geq {n+1}}, \;\;\PP-\mbox{a.s.}.
\end{align*}
The second component, $\xi = \Braces{ \xi_i, i \geq 1 }$, is a sequence of $\mathbf{R}$-valued random variables with transition probabilities
\begin{align*}
    \Prob{ \xi_n \in C \middle\vert \sigma\Braces{ (\tau_i, \xi_i), i = 1, \ldots, n-1 } \vee \sigma\Braces{ \tau_n } }
    =
    \frac{ K^\lambda_{\tau_n}(\omega, C ) }{ K^\lambda_{\tau_n}(\omega, \mathbf{R} ) }.
\end{align*}
Note that $(N^\lambda,\xi)$ can also be characterised by its so-called \textit{lifting}; that is, the random measure $\eta^\lambda$ defined as
\begin{align*}
    &\eta^\lambda(B) = \sum_{i \geq 1} \delta_{(\tau_i,\xi_i)}(B),
    &
    & B \in \Bb([0,T]) \otimes \Rr,
\end{align*}
with predictable compensator $K^\lambda_t(\omega, \diff r) \diff t$.

\begin{remark}
In general, when we refer to $(N^\lambda,\xi)$ as an $l$-dimensional marked Poisson process, we refer to the family
\begin{align*}
    (N^\lambda,\xi) = \Braces{ (N^{\lambda(j)},\xi^{(j)}), j = 1, \ldots, l },
\end{align*}
where each $(N^{\lambda(j)},\xi^{(j)})$ is an independent Poisson process with an $\FF$-intensity kernel $K^{\lambda(j)}$ from $(\Omega \times [0,T], \Ff \otimes \Bb([0,T]))$ into $(\mathbf{R},\Rr)$. Similarly, we write $K^{\lambda}$ (resp. $\eta^\lambda$) for the vector with $j$-th entry $K^{\lambda(j)}$ (resp. $\eta^{\lambda(j)}$); see (\ref{KerK}). 
\end{remark}

With the notion of admissible noise established above, we define next an \emph{admissible set-up}. Notice that within our framework, the random environment $\mu$ should be understood as the empirical measure process for the finite population interacting system. 
 
\begin{definition}
\label{Def:Admissible-Set-Up}
A tuple $(\Omega,\Ff,\FF,\PP,X_0,\mu,K^\lambda,(N^\lambda,\xi),W)$ (or simply $(X_0,\mu,K^\lambda,(N^\lambda,\xi),W)$) is said to be an \emph{admissible set-up} for the problem if
\begin{enumerate}
    \item
	$\mu$ is a random environment on $(\Omega, \Ff, \PP)$;  
    \item
	$X_0 \in L^2( \Omega, \Ff_0, \PP; \RR^d )$;
		
    \item
        $K^\lambda$ is a non-negative, $\FF$-predictable, locally integrable kernel from $(\Omega \times [0,T],\Ff \otimes \Bb([0,T])$ into the measurable space $(\mathbf{R},\Rr)$;
		
    \item
        $(N^\lambda,\xi)$ is an $l$-dimensional marked Poisson process, adapted to $\FF$, with $\FF$-stochastic intensity kernel $K^\lambda$ and lifting $\eta^\lambda$;

    \item
        A sequence  $W=\{W^i\}$ of independent $k$-dimensional $\FF$-Wiener processes, referred as an \emph{idiosyncratic noise};
		
    \item
        $(X_0, \mu, \eta^\lambda)$ is independent of $W$ (under $\PP$);
		
    \item
	Any $\FF^{(X_0, \mu, \eta^\lambda, W)}$-martingale is also an $\FF$-martingale (under $\PP$).
\end{enumerate}
\end{definition}

\begin{remark} 
\label{Remark:Intensity-Candidate}
\begin{enumerate}
    \item
        Condition 7 above is commonly known in the literature as \emph{$\Hh$-hypothesis}, meaning that $\FF^{(X_0, \mu, \eta^\lambda, W)}\hookrightarrow\FF$, or simply that ${(X_0, \mu, \eta^\lambda, W)}$ is \emph{compatible} with $\FF$; see \cite{aksamit_enlargement_2017}.

    \item 
        Since we are only considering kernels $K^\lambda$ that are obtained through an intensity candidate $\lambda$, we write $(X_0,\mu,\lambda,(N^\lambda,\xi),W)$ when referring to the admissible set-up $(X_0,\mu,K^\lambda,(N^\lambda,\xi),W)$.
    
    \item  
        Note that the   random environment $\mu$ involved in the solution of the SDEs (\ref{Eq:X^in})-(\ref{Eq:Underline_X}) is described in terms of the solution to be constructed for those SDEs.   
\end{enumerate}
\end{remark}

In order to describe properly  the last term  in (\ref{Eq:X^in})-(\ref{Eq:Underline_X}) corresponding to the integral with respect to the marked Poisson process, we need to define the following space:
\begin{align}
    \label{Eq:Frac-H-lambda}
    \mathfrak{H}^\lambda
    :=
    \Big\lbrace
        U : \Omega \times [0,T] \times \mathbf{R} \to \RR^{n \times l}
		\Big\vert\ &
		U \text{ is } \FF \otimes \Rr-\text{predictable, with } 
        \\ \nonumber
		&\Norm{ U_t(\omega,\cdot) }_{\lambda_t ( \omega )} < \infty,
        \\ \nonumber
		&\diff \PP \otimes \diff t\mathrm{-a.e. }(\omega, t) \in \Omega \times [0,T]
    \Big\rbrace,
\end{align}
where $\Norm{\ \cdot\ }_{\lambda_t (\omega)}^2 = \Parentheses{ \cdot, \cdot }_{\lambda_t (\omega)}$ is a random norm defined for matrix-valued functions, induced by the inner product
\begin{align}
    \label{Eq:lambda-inner-product}
    &\Parentheses{ u \cdot v }_{\lambda_t (\omega)} 
    := 
    \int_\mathbf{R} \Tr{ u(r) \diag{K^\lambda_t(\omega,\diff r)} v(r)^\top },
\end{align}
for $\diff \PP \otimes \diff t$-a.a. $(\omega, t) \in \Omega \times [0,T]$. This is due to the fact that for each $U \in \mathfrak{H}^\lambda$, if
\begin{align}
    \label{Eq:Finite-Predictable-Quadratic-variation-N^lambda}
    &\int_0^T \Norm{U_t(\cdot)}_{\lambda_t}^2 \diff t < \infty
    &
    &\PP-\mbox{a.s.},
\end{align}
then the integral process defined below is well-defined for $\diff \PP \otimes \diff t$-almost all $(\omega,t) \in \Omega \times [0,T]$:
\begin{align}
    \label{Eq:Integral-wrt-eta^lambda}
    \int_{[0,t] \times \mathbf{R} } U_s(\omega,r) {\eta}^\lambda &(\omega, \diff s, \diff r)
    :=
    \sum_{j=1}^l \sum_{i \geq 1} U^{(\cdot,j)}_{ \tau_i^{(j)} } (\omega,\xi_i^{(j)}(\omega)) \Ind{ \{\tau_i^{(j)} \leq t\} },
\end{align}
where $\tau^{(j)}$ (resp. $\xi^{(j)}$) is the sequence of jumping times (marks) of the $j$-member of $(N^\lambda,\xi)$, and $U^{(\cdot,j)}$ denotes the $j$-th column of $U$. Moreover, if the random variable in \eqref{Eq:Finite-Predictable-Quadratic-variation-N^lambda} is also integrable with repect to $\PP$, then the integrals with respect to the compensated martingale measure
\begin{align*}
    \tilde{\eta}^\lambda(\omega, \diff t, \diff r)
    :=
    {\eta}^\lambda(\omega, \diff t, \diff r)
    -
    K^{\lambda}_t(\omega, \diff r) \diff t
\end{align*}
are also  square-integrable martingales, see \cite[Thm. 5.1.33-34]{bremaud_point_2020}. 

We refer to \eqref{Eq:Integral-wrt-eta^lambda} as the \emph{integral of $U$ with respect to $(N^\lambda,\xi)$} --equivalently, the \emph{integral of $U$ with respect to $\eta^\lambda$}--. Following analogous notation  for well-posed integrals with respect to square-integrable martingales, we write
\begin{align*}
    \HH^2(\tilde{\eta}^\lambda)
    :=
    \Braces{
        U \in \mathfrak{H}^\lambda
        \  \middle\vert\ 
        \Esp{ \int_0^T \Norm{U_t(\cdot)}_{\lambda_t}^2 } < \infty
    }
\end{align*}
for the space of well-posed integrands with respect to $(N^\lambda,\xi)$, and 
\begin{align*}
    \Ll^2(\tilde{\eta}^\lambda)
    :=
    \Braces{
        \Parentheses{
            \int_{[0,t] \times \mathbf{R} } U_s(r) \tilde{\eta}^\lambda (\diff s, \diff r)
        }_{t \in [0,T]}
        \ \middle\vert\ 
        U \in \HH^2(\tilde{\eta}^\lambda)
    }
\end{align*}
for the space of martingales  expressed as compensated integrals with respect to $(N^\lambda,\xi)$.

\section{Existence of strong solutions}
\label{ExistenceStrongSol}
A key concept for the study of the  propagation of chaos is the notion of \emph{exchangeability}, which is used to  verify the convergence of a symmetric family of SDEs to a McKean-Vlasov system.
\begin{definition}
A family $X_0^I = \{ X^i_0, i \in I \}$ of random variables is said to be \emph{exchangeable} when the law of $X_0^I$ is invariant under any permutation of a finite number of indexes $i \in I$. 
\end{definition}

Observe that in the case of independent, identically distributed random variables, the exchangeability follows  immediately. The hypotheses about the model are presented below, introducing first those related with the probability space that support the noise processes,  and then presenting  the assumptions about the coefficients of the SDEs.
\begin{assumptions}
\label{Assumptions-on-the-Space}
\noindent
\begin{enumerate}
    \item \textbf{Probability space.} 
        There exists a complete probability space $(\Omega,\Ff,\PP)$ equipped with a general filtration $\FF=\{\Ff_t\}$, such that $l$ independent Poisson processes $N^{(i)}$ on $[0,T] \times \mathbf{R} \times \RR_+$ are supported ($(\mathbf{R},\Rr)$ being a measurable Polish space), with intensity measure $\diff t \otimes Q^{(i)}(\diff r) \otimes \diff s$. Moreover, this space supports a sequence $W=\{W^i\}$ of independent, $k$-dimensional Wiener processes.

    \item \textbf{Initial distribution.}
		There exists a measurable mapping 
		$X_0:(\Omega,\Ff) \to (\RR^d,\Bb(\RR^d))$ with law 
		$\nu_0 \in \Pp_2(\RR^d)$; i.e.
		\begin{align*}
            \Esp[^0]{ \Abs{ X_0 }^2 }=\int_{ \RR^d } \Abs{ x }^2 \nu_0( \diff x )
            <
            \infty.
		\end{align*}

    \item 
        \textbf{Intensity candidate.} There exists a measurable function $\lambda : \Omega \times [0,T] \times \Pp_2(\RR^d) \times \mathbf{R} \to \RR_+^l$, uniformly continuous in $[0,T]$ and Lipschitz continuous in $\Pp_2(\RR^d)$ (with respect to the $W_2$-distance), such that $\lambda$ is an admissible intensity candidate with respect to filtration $\FF^{\lambda,N} = \Braces{ \Ff_t^{\lambda,N}, 0 \leq t \leq T }$ defined as
        \begin{align*}
            \Ff_t^{\lambda,N}
            :=&
            \sigma\Braces{ \lambda_s(\nu,r): \forall (s,\nu,r) \in [0,t] \times \Pp_2(\RR^d) \times \mathbf{R} }
            \\
            &\vee
            \sigma\Braces{ N( B \times C \times D ), (B,C,D) \in \Bb([0,t]) \otimes \Rr \otimes \Bb(\RR_+)  }.
        \end{align*}
        
\end{enumerate}

\end{assumptions}

\begin{remark}
\label{Remark-on-Space}
The construction of the probability space $(\Omega,\Ff,\PP)$ in Assumption \ref {Assumptions-on-the-Space}.1 is assumed to be done by the product space $(\Omega^0 \times \Omega^1, \Ff^0 \otimes \Ff^1, \PP^0 \otimes \PP^1)$ as it is explained in \cite{hernandez-hernandez_coupled_2023}, with $\Omega^i,\;i=0,1$ being Polish spaces; this property turns out to be crucial in order to ensure that the regular conditional probabilities appearing in the next sections are properly defined, as well as  allowing us to adequately separate each of the noises affecting the SDEs under study. For the sake of completeness and readability, the details of the construction of this particular probability space and their consequences are  postponed to the Appendix \ref{Sec:Canonical-space}.
\end{remark}

\begin{assumptions}
\label{Assumptions-on-the-Coefficients}
Let
\begin{align*}
    &(b, \sigma):[0,T] \times \Pp_2(\RR^d) \times \RR^d \to \RR^d \times \RR^{d \times k},
    &
    &\gamma:[0,T] \times \Pp_2(\RR^d) \times \mathbf{R} \times \RR^d \to \RR^{d \times l}
\end{align*}
be measurable functions, uniformly continuous on $[0,T]$, such that the following conditions holds.
\begin{enumerate}
    
    \item \textbf{Lipschitz:} 
        There exists positive constants $K, K^0 > 0$, such that for any $x,x' \in \RR^d$, $\nu,\nu' \in \Pp_2(\RR^d)$, the following inequality holds for $\diff \PP \otimes \diff t$-a. a. $(\omega,t) \in \Omega \times [0,T]$:
		\begin{align*}
            &\Abs{ \delta b_t(\nu, x ) }^2
            +
            \Norm{ \delta \sigma_t(\nu, x )}^2_{\mathrm{Fr}}
            +
            \int_\mathbf{R} \Tr{ 
                \delta \gamma_t(\nu, r, x ) 
                \diag{\lambda_t(\nu,r)Q(\diff r)} 
                \delta \gamma_t(\nu, r, x )^\top 
            } 
            \\
		  &\quad\leq
		  K^0\ W_2(\nu,\nu')^2
		  +
		  K \Abs{ \delta x }^2.
		\end{align*}

    \item \textbf{Linear growth:} 
        There exists a positive constant $\beta > 0$ such that, for $(\omega,t)\in \Omega \times [0,T]$ fixed, 
        \begin{align*}
            &\Abs{ b_t( \nu, x ) }
            +
            \Norm{ \sigma_t( \nu, x ) }_{\mathrm{Fr}}
            +
            \Parentheses{ \int_\mathbf{R} \Tr{ 
                \gamma_t( \nu, r, x ) 
                \diag{\lambda_t(\nu,r)Q(\diff r)} 
                \gamma_t( \nu, r, x )^\top 
            } }^\frac{1}{2}
            \\
            &\quad\leq
		  \beta \Braces{ W_2(\nu,\delta_{\{0\}}) + \Abs{ x }}
        \end{align*}
        $\diff t \otimes \diff \PP$-a.e. for  $x \in \RR^d$ and  $\nu \in \Pp_2(\RR^d)$.
        
        \item {\bf Fourth moment:} 
            The coefficient $\gamma$ satisfies the following fourth-moment condition
            \begin{align}
                \label{Eq:gamma^*}
                &\max_{i \in \{1,\ldots,n\} } 
                    \int_\mathbf{R} \Abs{ \gamma^{(i,\cdot)}_t(\nu,x,r)}^4 Q^{(i)}( \diff r ) 
                \leq \gamma^*
                < \infty,
            \end{align}
            for each $(t,\nu,x) \in [0,T] \times \Pp_2(\RR^d) \times \RR^d$.      
\end{enumerate}
\end{assumptions}
Next, we establish the existence of solution in $(\Omega,\Ff,\PP)$ for both, the  \textit{finite-population system}  \eqref{Eq:X^in} and the McKean-Vlasov SDE \eqref{Eq:Underline_X}.

\begin{theorem}
\label{Thm:Strongexistence}
Suppose that Assumptions \ref{Assumptions-on-the-Space} and \ref{Assumptions-on-the-Coefficients} hold. Then,
\begin{enumerate}

    \item[(i)]
        there exists a unique strong solution process $X^{n} \in \mathcal{S}^2(\RR^{dn})$ to the finite-population system 
        \begin{align}
            \tag{\ref{Eq:X^in}}
            X^{i,n}_t
            =&
            X^{i,n}_0
            +
            \int_0^t b_s( \overline{ \mu }^{n}_{s-}, X^{i,n}_{s-} ) \diff s
            +
            \int_0^t \sigma_s( \overline{ \mu }^{n}_{s-}, X^{i,n}_{s-} ) \diff W_s^{i}
            \\ \nonumber
            &+
            \int_{[0,t] \times \mathbf{R}} 
                \gamma_s( \overline{ \mu }^{n}_{s-}, r, X^{i,n}_{s-})
            \tilde{\eta}^{\lambda,n}(\diff s, \diff r),\;\;\;t\leq T,
        \end{align}
        for  $i=1,\ldots,n$.
        
    \item[(ii)] 
        let $W^0$ be the first element of the sequence $W$, and let $X_0^0$ be an independent copy of $X_0$; then, there exists a unique strong solution process $X\in \mathcal{S}^2(\RR^d)$ to the conditional McKean-Vlasov (CMV) SDE with jumps 
        \begin{align}
            \label{Eq:CMV-X-1}
            X_t
            =&
            X_0^0
            + 
            \int_0^t b_s \Parentheses{ \Ll^1\Parentheses{ X }_{s-}, X_{s} } \diff s
            +
            \int_0^t \sigma_s \Parentheses{ \Ll^1\Parentheses{ X }_{s-}, X_{s} } \diff W^0_s
            \\ \nonumber
            &+
            \int_{[0,T] \times \mathbf{R} } \gamma_s \Parentheses{ \Ll^1\Parentheses{ X }_{s-}, r, X_{s-} } \tilde{\eta}^\lambda(\diff s, \diff r),
        \end{align}
        where $\eta^\lambda$ is the lifting of the marked Poisson process $(N^\lambda,\xi)$, with $N^\lambda$ a doubly stochastic Poisson process of intensity kernel
        \begin{align*}
            K^\lambda_t(\diff r) = \lambda_t(\Ll^1(X)_{t-},r) Q(\diff r),
        \end{align*}
        $\xi$ a sequence of marks on $(\mathbf{R},\Rr)$, and $\tilde{\eta}^\lambda$ the compensated martingale measure of $\eta^\lambda$.  
\end{enumerate}
\end{theorem}

Letting $X$ be the strong solution to the system \eqref{Eq:CMV-X-1} established in the above theorem, the pair $\Parentheses{\Ll^1(X), X}$ is defined as the \emph{solution process} for the conditional McKean-Vlasov SDE.
 
\begin{proof}
For part (i), observe that from the definition of Wasserstein distances, it follows that
\begin{align*}
    W_2(\overline{\mu}_x,\delta_{\{x_0\}})^2 \leq \frac{1}{n} \sum_{i=1}^n \Abs{x^i - x_0}^2,
\end{align*}
which together with Assumption \ref{Assumptions-on-the-Coefficients} imply that  \eqref{Eq:X^in} can be solved using  classical results from the theory of SDEs with random coefficients \cite[Thms. 14.21 \& 14.23]{jacod_calcul_2006}.

Now we proceed to prove part (ii).  Let  $W^0$ is defined on $(\Omega^1,\Ff^1,\FF^1,\PP^1)$ independently from $W$, see Remark \ref{Remark-on-Space}. For each $x \in \Ss^2( \RR^d )$,  define the mapping $\Phi : \Ss^2( \RR^d ) \to \Ss^2( \RR^d )$ as
\begin{align*}
    \Phi(x)_t
    :=&
    \chi
    +
    \int_0^t b_s \Parentheses{ \Ll^1\Parentheses{ x }_{s-}, x_{s} } \diff s
    +
    \int_0^t \sigma_s \Parentheses{ \Ll^1\Parentheses{ x }_{s-}, x_{s} } \diff W^0_s
    \\
    &+
    \int_{[0,t] \times \mathbf{R}} 
        \gamma_s \Parentheses{ \Ll^1\Parentheses{ x }_{s-}, r, x_{s} } 
    N\Parentheses{ \diff s, \diff r, [ 0, \lambda_s(x,r) ] },
\end{align*}
where we have use the slight abuse of notation $\lambda_s(x,r) = \lambda_s(\Ll^1\Parentheses{ x }_{s-},r)$, and $\chi \in L^2(\Ff_0;\RR^d)$.

Due to the Itô isometry with respect to $W^0$ and $\tilde{\eta}^\lambda$, as well as Burkholder-Davis-Gundy martingale inequality, the difference $\delta \Phi = \Phi(x) - \Phi(x')$ for any pair $x,x' \in \Ss^2(\RR^d)$ is bounded in norm $\Norm{\ \cdot\ }_*$ by
\begin{align*}
    \nonumber
    \Norm{ \delta \Phi }_*^2
    \leq&
    \Esp{ \int_0^T \Abs{ \delta b_t(x) }^2 + \Norm{ \delta \sigma_t(x) }^2_\mathrm{Fr} \diff t }
    \\
    &+
    \Esp{ \int_0^T \int_{ \mathbf{R} } \Tr{
        \delta \gamma_t(x,r) \diag{\lambda_t(x,r) Q(\diff r)} \delta \gamma_t(x,r)^\top
    } \diff t }
    \\
    &+
    \Esp{ \int_0^T \int_{ \mathbf{R} } \Tr{
        \gamma_t(x',r) \diag{ \Abs{\delta \lambda_t(x,r)} Q(\diff r)} \gamma_t(x',r)^\top
    } \diff t }.
\end{align*}
Let $\widehat{Q}^{(i)}(\diff r) =: \frac{Q^{(i)}(\diff r)}{ Q(\mathbf{R}) }$. Using Hölder's inequality,   last term of the RHS can be bounded above as follows:
\begin{align*}
    &\Esp{ \int_0^T \int_{ \mathbf{R} } \Tr{
        \gamma_t(x',r) \diag{ \Abs{\delta \lambda_t(x,r)} Q(\diff r)} \gamma_t(x',r)^\top
    } \diff t }
    \\
    &=
    \sum_{i=1}^d Q^{(i)}(\mathbf{R}) \Esp{ \int_0^T \int_{\mathbf{R}} 
        \Abs{ \gamma^{(i,\cdot)}_t(x',r) }^2 \Abs{ \delta \lambda^{(i)}(t,x,r) } 
    \widehat{Q}^{(i)}(\diff r) \diff t }
    \\
    &\leq
    \sum_{i=1}^d Q^{(i)}(\mathbf{R}) \EE \Bigg[  \int_0^T 
        \Parentheses{ \int_{\mathbf{R}} \Abs{ \gamma^{(i,\cdot)}_t(x',r) }^4 \widehat{Q}^{(i)}(\diff r)  }^{\frac{1}{2}}
        \\
        &\qquad\qquad\qquad\qquad\qquad\times 
        \Parentheses{ \int_{\mathbf{R}} \Abs{ \delta \lambda^{(i)}(t,x,r) }^2 \widehat{Q}^{(i)}(\diff r)  }^{\frac{1}{2}}
    \diff t \Bigg]
    \\
    &\leq
    d  \sqrt{K^* \gamma^*} \Abs{ Q(\mathbf{R}) }
        \Esp{ \int_0^T W_2( \Ll^1(x)_{t}, \Ll^1(x')_{t} ) \diff t }
\end{align*}

Then, by the Lipschitz conditions,
\begin{align*}
    \nonumber
    \Norm{ \delta \Phi }_*^2
    \leq&
    C_1 \int_0^T \Braces{
        \Esp{ W_2( \Ll^1(x)_t, \Ll^1(x')_t )^2 } 
        +
        \Esp{ \Abs{\delta x_t}^2 }
    } \diff t
    \\ \nonumber
    &+ C_2 \int_0^T \Esp{ W_2( \Ll^1(x)_{t}, \Ll^1(x')_{t} ) } \diff t
    \\ \nonumber 
    \leq&
    C \int_0^T \Braces{
        \Esp{ \Abs{ \delta x_t }^2 } + \Esp{ \Esp[^1]{ \Abs{ \delta x_t }^2 }^\frac{1}{2} }
    } \diff t
    \\ 
    \leq&
    C \int_0^T \Braces{ 
        \Esp{ \Abs{ \delta x_t }^2 } + \Esp{ \Abs{ \delta x_t }^2 }^\frac{1}{2} 
    } \diff t,
\end{align*}
where the second-to-last inequality is due to \eqref{Eq:W_2-Inequality}. Observe that
\begin{align}
    \label{Eq:delta-Phi-bound}
    \Norm{ \delta \Phi }_*^2
    \leq
    2C \max\Braces{ 
        \int_0^T \Esp{\Abs{ \delta x_t }^2} \diff t, 
        \Parentheses{ T \int_0^T \Esp{\Abs{ \delta x_t }^2} \diff t }^\frac{1}{2}
    }.
\end{align}

We proceed by cases. Suppose that the maximum in the above display is achieved in the second term, i.e.
\begin{align*}
    \Norm{ \delta \Phi }_*^2
    \leq
    2C \int_0^T \Esp{\Abs{ \delta x_t }^2} \diff t.
\end{align*}
Then, for any positive integer $k>1$, let $\Phi^k$ be the $k$-th composition of $\Phi$ with itself. In order to apply Banach's fixed point principle we need to verify that $\Phi^k$ is a contracting map for a sufficiently large $k$. Indeed, by iterating $k$ times we obtain that
\begin{align*}
    \Norm{ \delta \Phi^k }
    \leq&
    C \int_0^T \Esp{ \Abs{ \delta \Phi^{k-1}_t }^2 } \diff t
    \leq
    C \int_0^T \Esp{ \sup_{[0,t]} \Abs{ \delta \Phi^{k-1}_s }^2 } \diff t
    \\
    \leq&
    C^2 \int_0^T \int_0^t \Esp{ \Abs{ \delta \Phi^{k-2}_s }^2 } \diff s \diff t
    =
    C^2 \int_0^T \frac{T - t}{1} \Esp{ \Abs{ \delta \Phi^{k-2}_t }^2 } \diff t
    \\
    \leq&
    C^2 \int_0^T \frac{T-t}{1} 
    	\Esp{ \sup_{[0,t]} \Abs{ \delta \Phi^{k-2}_s }^2 } 
    \diff t
    \\
    \leq&
    C^3 \int_0^T \frac{T-t}{1} \int_0^t 
    	\Esp{ \Abs{ \delta \Phi^{k-3}_s }^2 } 
    \diff s \diff t
    =
    C^3 \int_0^T \frac{(T-t)^2}{2!} \Esp{ \Abs{ \delta \Phi^{k-3}_t }^2 } \diff t
    \\
    \leq& \ldots
    \\
    \leq&
    C^k \int_0^T \frac{(T-t)^{k-1}}{k!} \Esp{ \Abs{ \delta x_t }^2 } \diff t
    \leq
    \frac{(CT)^{k}}{k!} \Norm{ \delta x }_*.
\end{align*}

Now suppose that the maximum in the r.h.s. of \eqref{Eq:delta-Phi-bound} is achieved in the first term. Hence,
\begin{align*}
    \Norm{ \delta \Phi }_*^2
    \leq
    2C\sqrt{T} \Parentheses{\int_0^T \Esp{\Abs{ \delta x_t }^2} \diff t}^\frac{1}{2}.
\end{align*}
In this case, we have that
\begin{align*}
    &\Norm{ \delta \Phi }_*^2 \leq 2 C T \Norm{ \delta x }_*
\end{align*}
Let $\theta(x) = \frac{1}{2} \cosh{ \frac{x}{2CT}}$. Then, 
\begin{align*}
    \theta( \Vert \delta \Phi^k \Vert_* )
    &\leq
    \frac{1}{2} \Parentheses{ e^{ \Parentheses{ \frac{ \Vert \delta \Phi^{k} \Vert_* }{ 2CT } }^2 } }^\frac{1}{2}
    \leq
    \frac{1}{2}\Parentheses{ e^{ \frac{\Norm{ \delta \Phi^{k-1} }_*}{2CT} } }^\frac{1}{2}
    \\
    &\leq
    \frac{1}{\sqrt{2}} \Parentheses{ 
        \frac{
            e^{ \frac{\Norm{ \delta \Phi^{k-1} }_*}{2CT} } 
            +
            e^{ -\frac{\Norm{ \delta \Phi^{k-1} }_*}{2CT} } 
        }{
            2
        }
    }^\frac{1}{2}
    = 
    \theta( \Vert \delta \Phi^{k-1} \Vert_* )^\frac{1}{2}
    \\
    &\leq\ldots\leq
    \theta( \Vert \delta x \Vert_* )^\frac{1}{2^k}
\end{align*}
In other words, we have that $\Phi^k$ is a $\theta$-contraction mapping in the sense of  \cite[Definition 2.1]{ahmad_fixed_2017} with parameter $\frac{1}{2^k}$, and therefore admits a unique fixed point. Finally, by taking $\chi=X_0$, there exists a unique solution process $X$ to \eqref{Eq:CMV-X-1}, obtained as the fixed point of $\Phi^k$.

\end{proof}

\begin{remark}
Notice that the tuple
\begin{align}
    \label{Eq:Empirical-Admissible-Set-Up}
    \Parentheses{
        X_0,\overline{\mu}_{X^n},\lambda^n,(N^{\lambda,n},\xi^{\lambda,n}),(W^1,\ldots,W^n)
    }
\end{align} 
is not admissible in the sense of Definition \ref{Def:Admissible-Set-Up},
since Condition 6 is not satisfied. Nevertheless, this set-up will be complemented with an appropriate \textit{mean-field set-up} in which the conditional MacKean-Vlasov equation \eqref{Eq:CMV-X-1}  is solved, analyzing the convergence in the weak sense of both systems of SDEs.  
\end{remark}

\section{Conditional propagation of chaos}
\label{Sec:Propagation}

In this section the main theorem of this paper is presented, which consists in proving the \emph{propagation of chaos} from the finite population system \eqref{Eq:X^in} to the mean-field limit \eqref{Eq:CMV-X-1}. The way we proceed is by a \textit{synchronous coupling} between both systems (for a review into this subject, we refer to the monography by Chaintron \& Diez  \cite{chaintron_propagation_2022},\cite{chaintron_propagation_2022-1} and references there-in). Throughout this section we assume without explicit mention that Assumptions \ref{Assumptions-on-the-Space} and \ref{Assumptions-on-the-Coefficients} hold, and point out again that the probabilistic framework is the one described in the Appendix \ref{Sec:Canonical-space}.

\begin{theorem}
\label{Thm:Propagation-of-Chaos}
 For any finite collection of indexes $I \subset \NN$, the finite population system $X^{I,n}$ selected from \eqref{Eq:X^in} converges to the family of exchangeable McKean-Vlasov differential equations $X^I$ in the sense that
\begin{align}
    \label{Eq:(McK-V)Propagation-of-Chaos}
    \lim_{n \to \infty} \Braces{ 
        \Norm{ X^{I,n} - X^I }_*^2
        +
        \max_{i \in I}\ \sup_{ t \in [0,T] } \Esp{ W_2\Parentheses{ 
            \overline{\mu}_{t}^{i,n},
            \Ll^1\Parentheses{ X^i }_t
        }^2 }
    }
    =
    0.
\end{align}

\end{theorem}

\begin{remark}
Observe that due to the Wasserstein inequality \eqref{Eq:W_2-Inequality}, the first term in the LHS of the above display   implies that
\begin{align*}
    \lim_{n \to \infty} \sup_{[0,T]} W_2 \Parentheses{ \Ll(X^{i,n})_t, \Ll(X^i)_t } = 0
\end{align*}
for each $i \in I$; thus, from the relation between Wasserstein distances and weak convergence \cite[Thm. 5.11]{santambrogio_optimal_2015}, we conclude that $X^{I,n} \to X^I$ weakly when $n \to \infty$.
\end{remark}

Before proceeding to the proof of the above theorem, observe that the admissible set-up described in \eqref{Eq:(McK-V)Admissible-Set-Up} has the particular characteristic that the random environment is $\FF^0$-adapted; see (\ref{Eq:FF^0}). This property is useful  to prove the exchangeability of conditional McKean-Vlasov equations.

\begin{lemma}
\label{Lemma:McK-V-Exchangeability}
Let $I \subset \NN$ be a finite collection of indexes.  Then, the family of processes $X^I = \Braces{ X^i, i \in I }$
is conditionally independent, identically distributed given $\FF^0$.
Here,  $X^i$ is the solution process to \eqref{Eq:CMV-X-1} defined on the admissible set-up 
\begin{align*}
    \Parentheses{ \Omega, \Ff, \FF^0 \vee \FF^{W^i}, \PP, X_0^i, \Ll^1(X^i), \lambda^i, (N^{\lambda^i},\xi), W^i },
\end{align*}
with $\lambda^i:=\lambda(\cdot,\Ll^1(X^i),\cdot)$.  
\end{lemma}
\begin{proof}
For each $i \in I$, let $\mu^i:=\Ll^1(X^i)$ and $\FF^i = \FF^0 \vee \FF^{W^i}$. Notice that those coefficients  automatically satisfy the measurability and Lipschitz conditions, and  hence Theorem \ref{Thm:Strongexistence} yields the existence of  a strong solution $X$ to the MaKean-Vlasov equation \eqref{Eq:CMV-X-1}. Moreover, from Proposition \ref{Prop:Yamada-Watanabe}, there exists a mapping $\boldsymbol{\Phi}$, depending only on the law of the set-up, such that
\begin{align*}
    &\Prob{ X^i = \boldsymbol{\Phi}\Parentheses{ X_0^i, \mu^i, \lambda^i, (N^{ \lambda^i },\xi^i), W^i} } = 1,
    &
    &\forall i \in I.
\end{align*}

For any pair $(i,j) \in I \times I$, $i \neq j$, let $\underline{X}^{(i,j)}$ be the solution process to the auxiliary SDE
\begin{align*}
    \underline{X}_t^{(i,j)}
    =&
    X_0^i
    +
    \int_0^t b_s \Parentheses{ \mu^i_{s-}, \underline{X}^{(i,j)}_{s-} } \diff s
    +
    \int_0^t \sigma_s \Parentheses{ \mu^i_{s-}, \underline{X}^{(i,j)}_{s-} } \diff W_s^{j}
    \\
    &+
    \int_{[0,t] \times \mathbf{R} }\gamma_s \Parentheses{ \mu^i_{s-}, r, \underline{X}^{(i,j)}_{s-}} 
        \tilde{\eta}^{\lambda^i}(\diff s, \diff r),
    &
    &\forall t \in [0,T],
\end{align*}
The existence of a unique strong solution to this SDE follows from \cite[Thm. 3.3]{hernandez-hernandez_coupled_2023}, since $\mu^i$ is $\FF^0$-adapted, and the  measurability and Lipschitz conditions for the coefficients are satisfied. Therefore, by the uniqueness in law
\begin{align*}
    &\Prob{ 
        \underline{X}^{(i,j)} 
        = 
        \boldsymbol{\Phi}\Parentheses{ X_0^i, \mu^i, \lambda^i, (N^{\lambda^i},\xi^i), W^j }
    } = 1,
    &
    &\forall i \neq j.
\end{align*}

Since $W^i$ and $W^j$ are identically distributed on $(\Omega^1,\Ff^1,\PP^1)$, we obtain that for $\PP^0$-almost all $\omega^0 \in \Omega^0$ and for all $t \in [0,T]$,
\begin{align*}
    \mu^i_t(\omega^0)
    =
    \Ll^1\Parentheses{ X^i }_t(\omega^0)
    =
    \Ll^1\Parentheses{ \underline{X}^{(i,j)} }_t(\omega^0)
    =:
    \underline{\mu}^{(i,j)}_t(\omega^0),
\end{align*}
which in turn yields the equality between intensity candidates
\begin{align*}
    &\lambda^i = \lambda(\cdot,\mu^i,\cdot) = \lambda(\cdot,\underline{\mu}^{(i,j)},\cdot) 
        =: \underline{\lambda}^{(i,j)},
    &
    &\PP^0\mbox{-a.s.}
\end{align*}
Since the initial conditions are independent, identically distributed, they are also exchangeable, and hence 
\begin{align*}
    &\Prob{ 
        \underline{X}^{(i,j)} 
        = 
        \boldsymbol{\Phi}\Parentheses{ 
            X_0^i, \underline{\mu}^{(i,j)}, \underline{\lambda}^{(i,j)}, 
            (N^{ \underline{\lambda}^{(i,j)} },\xi^{(i,j)}), W^j 
        }
    }
    \\
    &=
    \Prob{ 
        \underline{X}^{(i,j)} 
        = 
        \boldsymbol{\Phi}\Parentheses{ 
            X_0^j, \underline{\mu}^{(i,j)}, \underline{\lambda}^{(i,j)}, 
            (N^{ \underline{\lambda}^{(i,j)} },\xi^{(i,j)}), W^j 
        }
    }
    = 
    1;
\end{align*}
in other words, $\underline{X}^{(i,j)}$ and $X^j$ share the same distribution, thereby proving that
\begin{align*}
    \Ll^1\Parentheses{ X^i }_t(\omega^0)
    =
    \Ll^1\Parentheses{ X^j }_t(\omega^0),
    &
    &\forall t \in [0,T].
\end{align*}
for $\PP^0$-a.a. $\omega^0 \in \Omega^0$.

Finally,  conditional independence of the elements of $X^ I$ is verified. Since  the underlying spaces are Polish, there exists a regular conditional probability $\mathrm{P}:\Omega^0 \times \Ff \to [0,1]$; that is, $\mathrm{P}$ is a transition kernel from $(\Omega^0,\Ff^0)$ into $(\Omega^1,\Ff^1)$ such that
\begin{align*}
    &\Prob{ (A^0 \times \Omega^1) \cap A } = \int_{A^0} \mathrm{P}(\omega^0,A) \diff \PP^0( \omega^0 )
\end{align*}
for all $A^0 \in \Ff^0$ and all $A \in \Ff$. However, since $\PP = \PP^0 \otimes \PP^1$,
\begin{align*}
    \mathrm{P}(\omega^0,A) 
    = 
    \Esp[^1]{ \Ind{A}(\omega^0,\cdot) } 
    = 
    \int_{\Omega^1} \Ind{A}(\omega^0,\omega^1) \diff \PP^1(\omega^1).
\end{align*}
Let $D^i, D^j \in \Dd([0,T];\RR^d)$. Then, from the definition of conditional expectation,
\begin{align*}
    &\int_{A^0 \times \Omega^1} \Esp{ 
        \Ind{ \Braces{ X^i \in D^i, X^j \in D^j } } 
        \middle\vert 
        \Ff^0 \otimes \{ \emptyset, \Omega^1 \} 
    }(\omega^0,\omega^1) \diff \PP(\omega^0,\omega^1)
    \\
    &=
    \int_{A^0 \times \Omega^1} \Ind{ \Braces{ X^i \in D^i, X^j \in D^j } }(\omega^0,\omega^1) 
        \diff \PP(\omega^0,\omega^1)
    \\
    &=
    \int_{A^0} \mathrm{P}\Parentheses{ \omega^0, \Braces{ X^i \in D^i, X^j \in D^j } } 
        \diff \PP^0( \omega^0 )
    \\
    &=
    \int_{A^0} \Prob[^1]{ X^i(\omega^0, \cdot ) \in D^i, X^j(\omega^0, \cdot ) \in D^j } 
        \diff \PP^0( \omega^0 )
    \\
    &=
    \int_{A^0} \Prob[^1]{ X^i(\omega^0, \cdot ) \in D^i } \Prob[^1]{ X^j(\omega^0, \cdot ) \in D^j } 
        \diff \PP^0( \omega^0 )
    \\
    &=
    \int_{A^0} \mathrm{P}\Parentheses{ \omega^0, \Braces{ X^i \in D^i } } 
        \mathrm{P}\Parentheses{ \omega^0, \Braces{ X^j \in D^j } } \diff \PP^0( \omega^0 )
    \\
    &=
    \int_{A^0 \times \Omega^1} 
        \Esp{ \Ind{ \Braces{ X^i \in D^i } } \middle\vert \Ff^0 \otimes \{ \emptyset, \Omega^1 \} }
            (\omega^0,\omega^1) 
        \\
        &\qquad\qquad\quad\times
        \Esp{ \Ind{ \Braces{ X^j \in D^j } } \middle\vert \Ff^0 \otimes \{ \emptyset, \Omega^1 \} }
            (\omega^0,\omega^1) 
    \diff \PP(\omega^0,\omega^1)
\end{align*}
for all $A^0 \in \Ff^0$,  concluding the proof.
\end{proof}

Finally, we are  ready to present the proof of the main result.

\begin{proof}[Proof of Theorem \ref{Thm:Propagation-of-Chaos}]
Without loss of generality, assume $\# I \leq n$ for some $n \gg 1$. From the proof for strong existence (Theorem \ref{Thm:Strongexistence}), there exists positive constants $C_1,C_2$ such that for all $i \in I$ and all $t \in [0,T]$,
\begin{align}
    \label{Eq:X^in-X^i}
    \Esp{ \sup_{[0,t]} \Abs{ X^{i,n}_s - X^i_s }^2 }
    \leq &
    C_1 \int_0^t \Braces{
        \Esp{ W_2( \overline{\mu}^{n}_s, \mu^{i}_s )^2 } 
        +
        \Esp{ \Abs{ X^{i,n}_s - X^{i}_s }^2 }
    } \diff s
    \\ \nonumber
    &+ 
    C_2 \int_0^t \Esp{ W_2( \overline{\mu}^n_s, \mu^{i}_s )^2 }^\frac{1}{2}  \diff s.
\end{align}
For all $j \in \{ 1, \ldots, n \} \setminus I$, let $X^j$ be a strong solution to the CMV \eqref{Eq:CMV-X-1} with respect to the Brwonian motion $W^j$ and write $\overline{\mu} := n^{-1} \sum_{i=1}^n \delta_{ X^i }$. Then, for $\omega^0$ fixed, 
\begin{align*}
    \Esp[^1]{ W_2( \overline{\mu}^n_t, \mu^{i}_t )^2 }
    &\leq
    \Esp[^1]{ W_2( \overline{\mu}^n_t, \overline{\mu}_{t} )^2 } + \Esp[^1]{ W_2( \overline{\mu}_{t}, \mu^{i}_t )^2 }
    \\
    &\leq
    \frac{1}{n} \sum_{i=1}^n \Esp[^1]{ \Abs{ X^{i,n}_t - X^i_t }^2 }
    +
    \Esp[^1]{ W_2( \overline{\mu}_{t}, \mu^{i}_t )^2 }
    \\
    &=
    \Esp[^1]{ \Abs{ X^{i,n}_t - X^i_t }^2 }
    +
    \Esp[^1]{ W_2( \overline{\mu}_{t}, \mu^{i}_t )^2 },
\end{align*}
where the last equality is due to the conditional exchangeability of $\Braces{ X^i, 1 \leq i \leq n }$ (Lemma \ref{Lemma:McK-V-Exchangeability}), and thus, the conditional exchangeability of $\{(X^{i,n},X^i), 1 \leq i \leq n \}$. Moreover, since $\Braces{ X^i, 1 \leq i \leq n }$ are conditionally independent, identically distributed given $\Ff^0$, from Glivenko-Cantelli theorem \cite[Thm. 5.8]{carmona_probabilistic_2018} it follows that for all $i \in I$
\begin{align}
    \label{Eq:Glivenko-Cantelli}
    &\Prob[^0]{ 
        \lim_{n \to \infty} \Esp[^1]{
            W_2\Parentheses{ \overline{\mu}_{t}, \Ll^1\Parentheses{ X^i }_t }^2
        } 
        = 0
    } = 1,
    &
    & \forall t \in [0,T].
\end{align}
However, observe that 
\begin{align}
    \nonumber
    \Esp[^1]{ W_2\Parentheses{ \overline{\mu}_{t}, \Ll^1\Parentheses{ X^i }_t }^2 }
    &\leq 
    2 \Esp[^1]{ 
        W_2 \Parentheses{ \overline{\mu}_{t}, \delta_0 }^2
        +
        W_2 \Parentheses{ \delta_0, \Ll^1\Parentheses{ X^i }_t }^2
    }
    \\ \label{Eq:(McK-V)Bound-difference-empiric-and-law}
    &\leq
    \frac{2}{n} \sum_{j=1}^n \Esp[^1]{ \Abs{ X^j_t }^2 }
    +
    2 \Esp[^1]{ \Abs{ X^i_t }^2 }
    =
    4 \Esp[^1]{ \Abs{ X^i_t }^2 }.
\end{align}
Given that $\Esp{ \Abs{ X^i_t }^2 } < \infty$, by Dominated Convergence we deduce that
\begin{align}
    \label{Eq:Pointwise-convergence-of-W_2}
    \lim_{n \to \infty} 
        \Esp{ W_2\Parentheses{ \overline{\mu}_{t}, \Ll^1\Parentheses{ X^i }_t }^2 } 
    &= 0,
    &
    &\forall t \in [0,T].
\end{align}

In order to verify that the convergence is uniform, let $s,t \in [0,T]$. From Cauchy-Schwartz and  {\eqref{Eq:(McK-V)Bound-difference-empiric-and-law}\ } we get
\begin{align*}
    &\Abs{ 
        \Esp{ W_2\Parentheses{ \overline{\mu}_{t}, \Ll^1\Parentheses{ X^i }_t }^2 }
        -
        \Esp{ W_2\Parentheses{ \overline{\mu}_{s}, \Ll^1\Parentheses{ X^i }_s }^2 }
    }^2
    \\
    &\leq
    \Esp{ \Braces{
        W_2\Parentheses{ \overline{\mu}_{t}, \Ll^1\Parentheses{ X^i }_t }
        +
        W_2\Parentheses{ \overline{\mu}_{s}, \Ll^1\Parentheses{ X^i }_s }
    }^2 }
    \\
    &\qquad\times
    \Esp{ \Braces{
        W_2\Parentheses{ \overline{\mu}_{t}, \Ll^1\Parentheses{ X^i }_t }
        -
        W_2\Parentheses{ \overline{\mu}_{s}, \Ll^1\Parentheses{ X^i }_s }
    }^2 }
    \\ 
    &\leq
    16 \Norm{ X^i }_*^2 \times
    \Esp{ \Braces{
        W_2\Parentheses{ \overline{\mu}_{t}, \Ll^1\Parentheses{ X^i }_t }
        -
        W_2\Parentheses{ \overline{\mu}_{s}, \Ll^1\Parentheses{ X^i }_s }
    }^2 }
    \\
    &\leq
    32 \Norm{ X }_*^2 \times \Braces{ 
        \Esp{ W_2\Parentheses{ \overline{\mu}_{t}, \overline{\mu}_{s} }^2 }
        +
        \Esp{ W_2\Parentheses{ \Ll^1\Parentheses{ X^i }_t, \Ll^1\Parentheses{ X^i }_s }^2
        }
    }
    \\
    &\leq
    64 \Norm{ X }_*^2 \times \Esp{ \Abs{ X^i_t - X^i_s }^2 },
\end{align*}
where the last inequality is due to \eqref{Eq:W_2-Inequality} and the conditional exchangeability of $\Braces{X^i, 1 \leq i \leq n}$. Then, due to the square-integrability of the \cadlag\ processes $X^i$ and $\Ll^1\Parentheses{ X^i }$, and the boundedness of the coefficients on bounded sets of $[0,T] \times \Pp_2(\RR^d) \times \RR^d$, there exists a positive constant $C'$ such that
\begin{align*}
    \Abs{ 
        \Esp{ W_2\Parentheses{ \overline{\mu}_{t}, \Ll^1\Parentheses{ X^i }_t }^2 }
        -
        \Esp{ W_2\Parentheses{ \overline{\mu}_{s}, \Ll^1\Parentheses{ X^i }_s }^2 }
    }^2
    \leq 
    C' \Abs{t - s}.
\end{align*}
Therefore, by equicontinuity, the convergence in \eqref{Eq:Pointwise-convergence-of-W_2} is uniform on $t$:
\begin{align}
    \label{Eq:Uniform-convergence-of-W_2-coupling}
    \lim_{n \to \infty}\ \sup_{[0,T]}
        \Esp{ W_2\Parentheses{ \overline{\mu}_{t}, \Ll^1\Parentheses{ X^i }_t }^2 } 
    &= 0,
    &
    &\forall i \in I.
\end{align}

Going back to \eqref{Eq:X^in-X^i}, let
\begin{align*}
    &b(t) := 1,
    &
    &k(s,t):= \max\Braces{C_1,C_2} \equiv k,
    &
    &g(t) := t + \sqrt{t}.
\end{align*}
Rewriting \eqref{Eq:X^in-X^i} in terms of $b,k$ and $g$, it follows that
\begin{align*}
    \Esp{ \sup_{[0,t]} \Abs{ X^{i,n}_s - X^i_s }^2 }
    \leq
    a(t)
    +
    b(t) \int_0^t k(t,s) g \Parentheses{ 
        \Esp{ \Abs{ X^{i,n}_s - X^{i}_s }^2 }
    } \diff s,
\end{align*}
where
\begin{align*}
    &a(t) := k t \sup_{s \in [0,t]} g\Parentheses{\Esp{ W_2\Parentheses{ \overline{\mu}_{s}, \mu^i_s }^2 }}.
\end{align*}
Applying the generalized non-linear version of Gronwall's inequality from \cite{butler_generalization_1971} (see also \cite[Ch. XII]{mitrinovic_inequalities_1991}) we get
\begin{align}
    \label{Eq:Gronwall-estimate}
    \Esp{ \sup_{[0,t]} \Abs{ X^{i,n}_s - X^i_s }^2 }
    &\leq
    G^{-1}\Parentheses{
        G(a(t))
        +
        k t
    },
\end{align}
where $G$ is the primitive integral 
\begin{align*}
    G(t) 
    = 
    \int^t_\epsilon \frac{\diff s}{g(s)}
    =
    2 \ln( \sqrt{t} + 1 ) - 2 \ln( \sqrt{\epsilon} + 1 ),
\end{align*}
and thus,
\begin{align*}
    G^{-1}(t) 
    = 
    \Parentheses{ \exp\Braces{ \frac{t}{2} + \ln( \sqrt{\epsilon} + 1 ) } - 1}^2
    =
    \Parentheses{ ( \sqrt{\epsilon} + 1 )e^\frac{t}{2} - 1}^2.
\end{align*}
However, since $g$ is positive for all $t > 0$, $G$ --and therefore $G^{-1}$-- is increasing in $[\epsilon,\infty)$ for all $\epsilon >0$. Furthermore, we have the estimates
\begin{align*}
    &G(t) \leq \int^t_\epsilon \frac{\diff s}{ s } = \ln( t ) - \ln( \epsilon )
\end{align*}
and
\begin{align*}
    &G^{-1}(t) 
    \leq 
    \Parentheses{ 2 \max\Braces{\sqrt{\epsilon}, 1} e^\frac{t}{2} }^2
    =
    4 \max\Braces{ \epsilon, 1} e^t.
\end{align*}

Going back to the estimate from \eqref{Eq:Gronwall-estimate}, observe that for $0 < \epsilon \ll 1$,
\begin{align*}
    &\Esp{ \sup_{[0,t]} \Abs{ X^{i,n}_s - X^i_s }^2 }
    \leq
    4 e^{kt} e^{G(a(t))}
    \leq
    \frac{ 4 e^{kT} }{\epsilon} a(T);
\end{align*}
then,
\begin{align*}
    \Esp{ \sup_{[0,t]} \Abs{ X^{i,n}_s - X^i_s }^2 }
    &\leq
    \inf_{\epsilon \in (0,1] } \Braces{ \frac{ 4 kT e^{kT} }{\epsilon} }
        \sup_{t \in [0,T]} g\Parentheses{\Esp{ W_2\Parentheses{ \overline{\mu}_{t}, \mu^i_t }^2 }}
    \\
    &=
    4 kT e^{kT} \sup_{t \in [0,T]} g\Parentheses{\Esp{ W_2\Parentheses{ \overline{\mu}_{t}, \mu^i_t }^2 }}
    = o(1),
\end{align*}
due to continuity of $g$ and the convergence of the coupling \eqref{Eq:Uniform-convergence-of-W_2-coupling}. In consequence
\begin{align}
    \label{Eq:Norm-star-convergence}
    &\lim_{ n \to \infty } 
    	\Esp{ \sup_{[0,t]} \Abs{ X^{i,n}_s - X^i_s }^2 } = 0,
    &
    &\forall t \in [0,T], \forall i \in I,
\end{align}
obtaining the first term in the asymptotic limit of \eqref{Eq:(McK-V)Propagation-of-Chaos}. Lastly, for the second term, note that
\begin{align*}
    \sup_{[0,T]} \Esp{ 
        W_2\Parentheses{ \overline{\mu}^{n}_t, \mu^i_t }^2 
    }
    \leq
    \underbrace{
        2 \sup_{[0,T]} \Esp{ \frac{1}{n} \sum_{i=1}^n \Abs{ X^{i,n}_t - X^{i}_t }^2 }
    }_{\to 0\ \eqref{Eq:Norm-star-convergence}}
    +
    \underbrace{
        2 \sup_{[0,T]} \Esp{ W_2\Parentheses{ \overline{\mu}_{t}, \mu^i_t }^2 }
    }_{\to 0\ \eqref{Eq:Uniform-convergence-of-W_2-coupling}},
\end{align*}
completing the proof.
\end{proof}

\begin{remark}
Notice that  the convergence of the empirical law $\overline{\mu}^n$ to the conditional law, instead of convergence to a  deterministic limit as in the classical case (the usual law of large numbers), is due to  the presence of a common noise in the system. Since the set-ups are only \textit{conditionally} independent and identically distributed (i.e. exchangeable), the use of Glivenko-Cantelli theorem in \eqref{Eq:Glivenko-Cantelli} guarantees the convergence to an $\Ff^0$-measurable distribution $\Ll^1(X)$; heuristically, this can be seen as  a deterministic limit law given that the common noise is known. Details on these claims can be found in  Appendix \ref{Sec:Canonical-space}.
\end{remark}

\section{Regime switching diffusions: An example}
\label{Sec:(Example)Regime-switching-diffusions}

Let $\Ee \subset \RR^n$ be a finite set of states and, for each $\nu \in \Pp_2(\RR^d)$, let $Q(\nu)=\Parentheses{ Q^{(i,j)}(\nu); i,j \in \Ee }$ be the transition semigroup of a continuous-time Markov chain; i.e. $Q$ is a conservative, irreducible matrix such that each entrance $Q^{(i,j)}$ is a mapping on $\Pp_2(\RR^d)$, which is Lipschitz with respect to the metric  $W_2$. Furthermore, assume that the mapping 
\begin{align*}
	\Ee \times \Pp_2(\RR^d)
	\ni
	(i,\nu)
	\mapsto
	Q^{(i)}(\nu)
	:=
	\sum_{j \in \Ee} Q^{(i,j)}(\nu)
\end{align*}
is uniformly bounded, i.e. there exists a positive constant $H_0$ such that 
\begin{align*}
	H_0 := \sup_{ \nu \in \Pp_2(\RR^d) } \sup_{i \in \Ee} Q^{(i)}(\nu) < \infty.
\end{align*}
Assume that the set  $\Ee \times \Ee$ is ordered according to the lexicographic ordering \cite{yin_hybrid_2010}. For all $i,j \in \Ee$ with $i \neq j$ and $\nu \in \Pp_2(\RR^d)$, let $\Gamma^{(i,j)}$ be the consecutive   left-closed, right-open intervals of $\RR_+$, each one having length $Q^{(i,j)}(\nu)$, and $\Gamma^{(i,i)} = \emptyset$:
\begin{align*}
	\Gamma^{(i,j)}(\nu)
	:&=
	\left[
		\underline{Q}^{( i, j )}(\nu)
		,
		\underline{Q}^{( i, j )}(\nu)
		+
		Q^{( i, j )}(\nu)
	\right),
\end{align*}
where
\begin{align*}
	\underline{Q}^{( i, j )}(\nu)
	:&=
	\underset{ \mathfrak{i} \neq \mathfrak{j} }{
		\sum_{ \mathfrak{i} \leq i }
		\sum_{ \mathfrak{j} \leq \underline{j}( \mathfrak{i}; i ) }
	}
	Q^{( \mathfrak{i}, \mathfrak{j} )}(\nu),
	&
	&\mbox{with}	
	&
	\underline{j}( \mathfrak{i}; i )
	&=
	\begin{cases}
		n,		& \mathfrak{i} < i,
		\\
		j-1,	& \mathfrak{i} = i.
	\end{cases}
\end{align*}

Given a random environment $\mu$ and an independent Poisson process $N$ of intensity measure $\diff s \otimes \diff r$, let $(N^Q,\xi^Q)$ be an admissible common noise defined (iteratively) by the lifting 
\begin{align*}
    \eta^Q(\diff s, \diff e)
    =&
    \sum_{j \in \Ee} N\Parentheses{
        \diff s, \Gamma^{ (i,j) }(\mu_{s-}) 
    } \otimes \delta_{\{j\}}(\diff e),
\end{align*}
on the event $\{ \xi^Q_{N^Q_{s-}} = i \}$, for each $i \in \Ee$. That is, $(N^Q,\xi^Q)$ admits an intensity kernel given by
\begin{align*}
    K^Q(s,\diff e)
    =
    \sum_{j \in \Ee} Q^{(i,j)}(\mu_{t-})\ \delta_{\{j\}}(\diff e)
    =
    Q^{( i, \cdot )}(\mu_{t-}) \delta_{\Ee}(\diff e).
\end{align*}
The process $Y^Q = \Braces{ Y_t, t \in [0,T] }$, defined as
\begin{align}
    \label{Eq:Y^Q}
    Y^Q_t
    :=
    \xi_{N^Q_t}^Q
    =
    \int_{ \{t\} \times \Ee } e \ \eta^Q(\diff s, \diff e),
\end{align}
is then referred to as the \textit{environment dependent regime}. That is, $Y^Q$ is a conditional Markov chain with transition matrix $Q$:
\begin{align*}
    \Prob{ Y^Q_{ t + \delta} = j \middle\vert Y^Q_t = i, \Ff^\mu_t  }
    =
    \begin{cases}
		Q^{(i,j)}(\mu_t)\delta + o(\delta),		& i \neq j,
		\\
		1 + Q^{(i,i)}(\mu_t)\delta + o(\delta),	& i=j.
    \end{cases}
\end{align*}

Let $b:\Pp_2(\RR^d) \times \RR^d \times \Ee \to \RR^d$ be a bounded function such that $(\nu,x)\mapsto b(\nu,x,e)$ satisfy the corresponding \textbf{Lipschitz} conditions for all $e$; let $\sigma$ be a real constant and $\gamma \equiv 0$. Then, $\Parentheses{\Omega,\Ff,\FF,\PP,\overline{X}^n_0,\overline{\mu}^n,K^n,(N^{n},\xi^n),W}$ forms an empirical admissible set-up, where $X^{\{1,\ldots,n\},n}$ is the solution process to
\begin{align}
    \label{Eq:(Example-II)X^n}
    X^{i,n}_t 
    &= 
    X_0^{i,n}
    +
    \int_0^t b \Parentheses{  \overline{\mu}^n_s, X^{i,n}_s, Y_s^n } \diff s
    +
    \sigma W^{i}_t,
    &
    &\forall i =1, \ldots,n,
\end{align}
where $Y^n$ is the regime determined (as in \eqref{Eq:Y^Q}) by the intensity kernel $K^n(t,\diff e)=Q(\overline{\mu}^n)\delta_\Ee(\diff e)$. In other words, $X^{\{1,\ldots,n\},n}$ above constitutes a state-dependent, regime-switching diffusion with transition matrix $Q( \overline{\mu}^n )$ \cite{yin_hybrid_2010}, \cite{shao_strong_2015}, \cite{shao_propagation_2022}.

Taking a set $I \subseteq \{ 1, \ldots, n \}$ fixed, from Theorem \ref{Thm:Propagation-of-Chaos} we know that $X^{I,n}$, converges to an exchangeable collection of conditional McKean-Vlasov regime-switching diffusions $X^I$ when $n \to \infty$, where each $X^i$ is a conditionally independent, identically distributed copy of
\begin{align}
    \label{Eq:(Example-II)Underline_X}
    X_t 
    &= 
    X_0 
    +
    \int_0^t b\Parentheses{ \Ll^1\Parentheses{ X }_s, X_s, Y_s }\diff s
    +
    \sigma W_t,
\end{align}
subjected to the regime $Y = \Braces{ Y_t := \xi_{N^{\lambda}_t}, t \in [0,T] }$, where $(N^{\lambda},\xi)$ is the admissible common Poissonian noise with intensity kernel defined on the event $\{ Y_{t-} = j \}$, for each states $j \in \Ee$ as
\begin{align}
    \label{Eq:(Example-II)Underline_lambda}
    K^{\lambda}(t,\diff e)
    &=
    Q^{( i, \cdot )}(\Ll^1(X)_{t-}) \delta_{\Ee}(\diff e).
\end{align}
Indeed, since each
\begin{align*}
    X^{i,n}_0 
    +
    \int_0^t b\Parentheses{ \overline{\mu}^n_s, X^{i,n}_s, j }\diff s
    +
    \sigma W^i_t
\end{align*}
converges weakly to
\begin{align*}
    X^{i}_0 
    +
    \int_0^t b\Parentheses{ \Ll^1\Parentheses{ X^i }_s, X^{i}_s, j }\diff s
    +
    \sigma W^i_t
\end{align*}
for all $j \in \Ee$, by introducing the intermediate coupling 
\begin{align*}
    \underline{X}^{i,n}_t 
    &= 
    X_0^{i,n}
    +
    \int_0^t b \Parentheses{  \overline{\mu}^n_s, X^{i,n}_s, Y_s } \diff s
    +
    \sigma W^{i}_t
\end{align*}
and conditioning with respect to the process $N$, the result is obtained.

Observe that these  results on the conditional propagation of chaos are consistent with those obtained by \cite{shao_propagation_2022}. Moreover, they also allow for a direct interaction between the state of the particle system and the regime itself.

\appendix
\section{A digression on the conditional law as a random environment}
\label{Sec:Canonical-space}

Let
\begin{align}
    \label{Eq:(Omega,Ff,PP)}
    (\Omega, \Ff, \PP) 
    := 
    \Parentheses{
        \Omega^0 \times \Omega^1, (\Ff^{\lambda,N}_T \otimes \Ff^{W}_T) \vee \Nn, \PP^0 \otimes \PP^1
    }
\end{align}
be described in the following way. Here $(\Omega^i,\Ff^i,\PP^i), i=0,1$, are probability spaces, with $\Omega^i$ being Polish spaces. Moreover,  $\lambda$ and $N$ are defined on $\Omega ^0$ and $W$ on $\Omega ^1$, respectively. We denote by  $\Nn$  the set $\PP^0\otimes\PP^1$-null sets on the product space. Similarly, the corresponding  filtrations are extended to the product space:
\begin{align}
    \label{Eq:FF^0}
    \FF^0 
    :&= 
    \Parentheses{ \FF^{\lambda,N} \otimes \Braces{ \emptyset, \Omega^1 } } \vee \Nn
    =
    \Braces{ ( \Ff_t^{\lambda,N} \otimes \Braces{ \emptyset, \Omega^1 } ) \vee \Nn=:\Ff^0_t, 0 \leq t \leq T },
    \\ \label{Eq:FF^1}
    \FF^1 
    :&= 
    \Parentheses{ \Braces{ \emptyset, \Omega^0 } \otimes \FF^{W} } \vee \Nn
    =
    \Braces{ ( \Braces{ \emptyset, \Omega^0 } \otimes \Ff_t^{W} ) \vee \Nn=:\Ff^1_t, 0 \leq t \leq T },
    \\ \label{Eq:FF}
    \FF :&= \FF^{0} \vee \FF^{1}.
\end{align}
Note that $\FF$, $\FF^0$ and $\FF^1$ satisfy the usual conditions, and $\FF^0$ and $\FF^1$ are independent. Furthermore, $\FF = \Parentheses{ \FF^{\lambda,N} \otimes \FF^{W} } \vee \Nn$.

For stochastic processes defined on each  coordinate space, we abuse slightly of notation and use the same symbols as their natural extension, i.e.
\begin{align}
    \label{Eq:Theta-extended}
    \Theta_t (\omega^0, \omega^1) :&= \Theta_{t} (\omega^0)
    &
    &\forall (t, \omega^1) \in [0,T] \times \Omega^1,
    \\ \nonumber
    &\Parentheses{\mbox{resp. }\Theta_{t}(\omega^1)},
    &
    &\Parentheses{\mbox{resp. }\forall (t, \omega^0) \in [0,T] \times \Omega^0}.
\end{align}
We also denote the conditional law of any stochastic process $\Theta$ in $\RR^d$ given the common noise, as the measured valued process $\Ll^1(\Theta) = \Braces{ \Ll^1(\Theta)_t, t \in [0,T] }$ given by
\begin{align}
    \label{Eq:Conditional-Ll^1}
    &\Ll^1(\Theta)_t(B):= \Prob{ \Theta_t \in B \middle\vert \Ff^0_T },
    &
    &\forall B \in \Bb(\RR^d).
\end{align}

In general, observe that within the current set-up, given a random variable $\Theta:(\Omega, \Ff)\to \Xx$, with $\Xx$  a Polish space, the mapping
\begin{align*}
    \Omega^1 
    \ni
    \omega^1 
    \mapsto
    \Theta( \omega^0, \omega^1 )
\end{align*}
is a random variable on $(\Omega^1,\Ff^1,\PP^1)$ for $\PP^0$-almost all $\omega^0 \in \Omega^0$. Define $\Nn^\Theta \subset \Omega^0$ as the complement of the set in which $\Theta( \omega^0,\ \cdot\ )$ is well-defined. Therefore, we can define $\tilde{\Ll}^1(\Theta):\Omega^0 \to \Pp(\Xx)$ as
\begin{align}
    \label{Eq:tilde-Ll^1}
    \omega^0 
    &\mapsto 
    \tilde{\Ll}^1(\Theta)
    =
    \begin{cases}
        \Ll( \Theta( \omega^0,\ \cdot\ ) ), & \omega^0 \not\in \Nn^\Theta,
        \\
        \delta_{\{0\}},                     & \omega^0 \in \Nn^\Theta.
    \end{cases}
\end{align}
Using \cite[Lemma 2.4]{carmona_probabilistic_2018-1}, this mapping is indeed a $\Pp(\Xx)$-valued random variable on $\Omega^0$, and coincides with our previous definition of $\Ll^1(\Theta)$ for random processes in \eqref{Eq:Conditional-Ll^1} $\PP^0$-a.s. thanks to the existence of a regular conditional probability measure. 

In particular, when $\Theta$ is a random variable taking values on the Polish space of trajectories $\mathrm{D}([0,T];\RR^d)$, then $\Ll^1(\Theta)$ is a $\Pp(\mathrm{D}([0,T];\RR^d))-$ valued random variable. Since we endowed $\Pp(\mathrm{D}([0,T];\RR^d))$ with the $\sigma$-algebra generated  from  measurable projections, we have that for each $t \in [0,T]$,
\begin{align*}
    &\Ll^1(\Theta)_t =\Ll^1(\Theta_t),
    &
    &\PP^0-\mbox{a.s.}.
\end{align*}
Furthermore, if $\Theta$ is an $\FF$-adapted      (respectively, square-integrable with continuous paths) stochastic process, then the corresponding $\Pp(\RR^d)$-valued process $\Ll^1(\Theta)$ is $\FF^{0}$-adapted (resp. square-integrable with continuous paths) \cite[Lemma 2.5]{carmona_probabilistic_2018-1}. These facts lead to an intrinsic relation between the discontinuities of both processes $\Theta$ and $\Ll^1(\Theta)$.

\begin{lemma}
\label{Lemma:Ll^1-is-cadlag-adapted}
Let $\Theta$ be an $\RR^d$-valued process with \cadlag\ trajectories. If $\Theta$ is square integrable and $\FF$-adapted, then there exists a version of $\Ll^1(\Theta)$ with \cadlag\ trajectories in $\Pp_2(\RR^d)$ (with respect to the $2$-Wasserstein distance),  $\FF^{0}$-adapted, and its jumping times match those of $\Theta$.
\end{lemma}

\begin{proof}
Observe that for all $t \in [0,T)$ and $\omega^0 \not\in \Nn^{\Theta}$, where $\Nn^{\Theta} \subset \Omega^0$ is a set of null measure,
\begin{align*}
    &\Theta_{t}(\omega^0,\omega^1) = \lim_{n \to \infty} \Theta_{t + \frac{1}{n}}(\omega^0,\omega^1),
    &
    &\forall \omega^1 \in \Omega^1.
\end{align*}
From our previous discussion, this implies that $\Parentheses{ \Ll^1(\Theta)_t }$ has $\PP^0$-a.s. right-continuous paths in $\Pp_2(\RR^d)$ with respect to the $2$-Wasserstein distance. Similarly, by defining 
\begin{align*}
    \Theta_{t-}(\omega^0,\omega^1)
    &:=
    \lim_{n \to \infty} \Theta_{ t - \frac{1}{n} }(\omega^0,\omega^1),
    &
    &\forall 
    (\omega^0,\omega^1) 
    \in 
    (\Omega^0 \setminus \Nn^{\tilde{\Theta}}) \times \Omega^1,
\end{align*}
which we know that exists $\PP$-a.s. for all $t \in (0,T]$, and setting 
\begin{align*}
    \omega^0 
    &\mapsto 
    \tilde{\Ll}^1(\Theta)_{t-}
    :=
    \begin{cases}
        \Ll( \Theta( \omega^0,\ \cdot\ )_{t-} ),		& \omega^0 \not\in \Nn^\Theta,
        \\
        \delta_{\Braces{0}},						& \omega^0 \in \Nn^\Theta,
    \end{cases}
\end{align*}
for each $t \in (0,T]$, we have obtained a version of $\Parentheses{ \Ll^1(\Theta)_t }$ which paths have left limits $\PP^0$-a.s. in $\Pp_2(\RR^d)$ with respect to the $2$-Wasserstein distance,  concluding the proof.

\end{proof}

Applying the previous result to the solution $X$ of the system \eqref{Eq:CMV-X-1}, we  conclude that $\Ll^1(X)$ is admissible with respect to $\FF$ as a random environment.

\begin{corollary}
\label{Lemma:McK-V-Admissible-Set-Up}
Under the same assumptions of Theorem \ref{Thm:Strongexistence},
\begin{align}
    \label{Eq:(McK-V)Admissible-Set-Up}
    \Parentheses{ X_0^0, \Ll^1(X), \lambda, (N^\lambda,\xi), W^0 }
\end{align}
forms an admissible set-up, in terms of Definition \ref{Def:Admissible-Set-Up}, where  $X$ is the solution process to the McKean-Vlasov equation \eqref{Eq:CMV-X-1}.
\end{corollary}

\begin{proof}
Observe that Conditions 2-7 from  Definition \ref{Def:Admissible-Set-Up} follows  directly from the construction of the probability space. In fact,  the compatibility condition is obtained  from the fact that we are considering the product space \eqref{Eq:(Omega,Ff,PP)}-\eqref{Eq:FF}, with the corresponding noises being independent. 
 
Hence, we only need  to prove that $\Ll^1(X)$ constitutes a random environment, i.e. a \cadlag,\ square-integrable, $\FF$-adapted measure-valued process. Let $\nu_0:[0,T] \to \RR^d$ be the constant process identically zero. Due to the definition of the $2$-Wasserstein distance, it follows that
\begin{align}
    \label{Eq:W_2-Inequality'}
    W_{2,d_\infty}\Parentheses{ \Ll^1\Parentheses{ X }, \nu_0 }^2
    \leq&
    \int_{ \mathrm{D}([0,T];\RR^d)^2 } 
        d_\infty(x,x')^2
    \Parentheses{ 
        \Ll^1\Parentheses{ X } \otimes \nu_0
    }( \diff x, \diff x' )
    \\ \nonumber
    =&
    \Esp[^1]{ d_\infty( X, 0 )^2 },\;\;\PP^0-a.s.
\end{align}
where $d_\infty$ is the distance induced by the supremum norm. As a consequence, Lemma \ref{Lemma:Ll^1-is-cadlag-adapted} and the square integrability of $X$ imply together that $\Ll^1(X)$ is a random environment in the probability space described in (\ref{Eq:(Omega,Ff,PP)}). 

\end{proof}

From the theory of stochastic integration recall that, under suitable conditions, strong solutions to SDEs can be regarded as functionals of the underlying noise processes (see, e.g. \cite[Thm. 2]{barczy_yamada-watanabe_2015} for SDEs with jumps). This is  known in the literature as the \textit{Yamada-Watanabe theorem} due to the original work done by these authors on SDEs driven by a Brownian motion \cite{yamada_uniqueness_1971}. For McKean-Vlasov type of equations this result is technically complicated, in particular when the SDE is conditional to a common noise filtration. However, the  construction of the probabilistic model presented in this work allow us  to work under the framework established by Carmona and Delarue  \cite[Sec. 1.2.3]{carmona_probabilistic_2018-1} with relatively minor changes to account for the jumps in the system. Namely, by considering a sufficiently rich input space, all the topological properties required for the corresponding version of the Yamada-Watanabe theorem to hold are satisfied. 

Let 
\begin{align}
    \label{Eq:Omega-input} 
    \Omega^{\mathrm{input}}
    :=&
    \RR^d \times \mathrm{D}( [0,T]; \Pp_2(\RR^d)) \times \Mm \Parentheses{ [0,T] \times \mathbf{R} }^l
        \\ \nonumber
        &\times \Mm_c^* \Parentheses{ [0,T] \times \mathbf{R} \times \RR_+ }^l
        \times \mathrm{C} \Parentheses{ [0,T]; \RR^k },
\end{align}
such that each marginal is endowed with its corresponding topology. Namely:
\begin{itemize}
    
    \item 
        $\RR^d$ with the usual Euclidean topology;
        
    \item
        $\mathrm{D}( [0,T]; \Pp_2(\RR^d) )$ with the Skorohod $J1$-topology with respect to the 2-Wasserstein distance $W_2$ \cite{dawson_measure-valued_1993};
        
    \item
        $\mathrm{C}( [0,T]; \RR^k )$ with the uniform topology \cite{billingsley_convergence_2013};
        
    \item
        $\Mm( [0,T] \times \mathbf{R} )$ and $\Mm_c^*( [0,T] \times \mathbf{R} \times \RR_+ )$ with the weak$^*$-topology \cite{daley_introduction_2008}.
        
\end{itemize}
Then, we define 
\begin{align}
    \label{Eq:Canon-space}
    (\boldsymbol{\Omega}, \boldsymbol{\Gg}, \mathbf{Q})
    =
    \Parentheses{ 
        \Omega^\mathrm{input} \times \mathrm{D}( [0,T]; \RR^d ), 
        \Bb(\Omega^\mathrm{input}) \otimes \Dd( [0,T]; \RR^d ),
        \mathbf{Q}
    }
\end{align}
as the measurable space with the product $\sigma$-algebra of the corresponding Borel $\sigma$-algebras and $\mathbf{Q}$ is the probability distribution induced by the random system
\begin{align}
    \label{Eq:Canon-Process}
    \mathfrak{C}
    :=&
    \Parentheses{X_0^0, \Ll^1(X), \int_{\cdot\times \cdot} K^\lambda_s(\diff r) \diff s, N, W^0, X }.
\end{align}
Furthermore, let $\mathbf{G}$ be the augmented filtration generated by
\begin{align*}
    \boldsymbol{\Gg}_t
    :=&
    \sigma\Braces{ X_0^0 }
    \vee
    \sigma\Braces{ 
    	\Ll^1(X)_s; 0 \leq s \leq t 
    }
    \\ \nonumber
    &\vee
    \sigma\Braces{ \int_{ B \times C }K^\lambda_s(\diff r) \diff s; B \in \Bb([0,t]), C \in \Rr  }
    \\ \nonumber
    &\vee
    \sigma\Braces{ 
        N( B \times D ); B \in \Bb([0,t]), D \in \Rr \otimes \Bb(\RR_+) 
    }
    \\ \nonumber
    &\vee
    \sigma\Braces{ W^0_s; 0 \leq s \leq t }
    \vee
    \sigma\Braces{ X_s; 0 \leq s \leq t },
\end{align*}
for all $t \in [0,T]$. The reason to choose this space, in particular the selection  of $\Mm \Parentheses{ [0,T] }^l$ for the coordinate related to the intensity process, is that we can gain some regularity by invoking the Lebesgue differentiation theorem. Indeed, given an admissible kernel $K^\lambda$ of the form $\lambda_t(r)Q(\diff r) \diff t$, we can recover a regularized version of the kernel by means of the transformation
\begin{align*}
    \tilde{K}^{\lambda}_t(C)
    :=
    \lim_{ n \to \infty } 
        n \int_{ \left( t - \frac{1}{ n }, t \right] \times C } \lambda_s(r) Q(\diff r)\diff s,
\end{align*}
if the limit exists, and $0$ in any other case. Thus, instead of   work directly with the intensity candidate $\lambda$ in the canonical space, we can use the measure $K^\lambda_t(\diff r) \diff t$, corresponding to the absolutely continuous compensator of the admissible noise.

\begin{proposition}
\label{Prop:Yamada-Watanabe}
The distribution law induced in $\Parentheses{ \boldsymbol{\Omega}, \boldsymbol{\Gg}, \mathbf{Q} }$ (see  \eqref{Eq:Canon-space}) by the random element $\mathfrak{C}$ defined in \eqref{Eq:Canon-Process},  depends only upon the law of the set-up \eqref{Eq:(McK-V)Admissible-Set-Up} 
\begin{align}
    \label{Eq:Input-Process}
    \Parentheses{X_0^0, \Ll^1(X), \int_{\cdot\times\cdot} K^\lambda_s(\diff r) \diff s, N, W^0 }
\end{align}
Furthermore, there exists a $\Bb(\Omega^\mathrm{input}) / \Dd( [0,T]; \RR^d ) $-measurable mapping 
\begin{align*}
    \boldsymbol{\Phi}:\Omega^\mathrm{input}\to\mathrm{D}([0,T];\RR^d),
\end{align*}
depending only on the distribution of the inputs \eqref{Eq:Input-Process}, such that the equality
\begin{align*}
    &X
    =
    \boldsymbol{\Phi}\Parentheses{ 
        X_0^0, \Ll^1(X), \int_{\cdot\times\cdot} K^\lambda_t(\diff r) \diff t, N, W^0 
    }
\end{align*}
holds $\PP$-a.s.
\end{proposition}

\begin{proof}
Let $\Big( \boldsymbol\chi, \mathfrak{M}, \boldsymbol\Lambda, \mathbf{n,w,x} \Big)$ be the canonical process on $\Parentheses{ \boldsymbol{\Omega}, \boldsymbol{\Gg}, \mathbf{G}, \mathbf{Q} }$ induced by $\mathfrak{C}$. Observe that, by construction, we can follow  exactly the same steps as in \cite[Theorem 1.33]{carmona_probabilistic_2018-1}, setting the backward components equal to zero. Therefore,   it is sufficient to verify that $\mathbf{n}$ and $\mathbf{w}$ are mutually independent on the canonical space. Indeed, from the Watanabe and L\'evy characterizations, respectively, it is  known that the former is a $\mathbf{G}$-Poisson process and the latter is a $\mathbf{G}$-Wiener process. Since both of them are defined on the same filtration where they are already independent, and hence this property is inherited \cite[Thm. II.6.3]{ikeda_stochastic_1989}.
\end{proof}

\begin{remark}
In general, notice that the previous proposition can be stated in terms of any strong solution process 
\begin{align*}
    X_t 
    =&
    X_0^0
    +
    \int_0^t b_s(\mu_{s-},X_{s-}) \diff s
    +
    \int_0^t \sigma_s(\mu_{s-},X_{s-}) \diff W_s^0
    \\
    &+
    \int_{[0,t] \times \mathbf{R}} \gamma_s(\mu_{s-},r,X_{s-}) \tilde{\eta}^\lambda( \diff s, \diff r)
\end{align*}
defined on an admissible set-up $\Parentheses{X_0^0, \mu, \lambda, (N^\lambda,\xi), W^0 }$.
\end{remark}




\providecommand{\bysame}{\leavevmode\hbox to3em{\hrulefill}\thinspace}
\providecommand{\MR}{\relax\ifhmode\unskip\space\fi MR }
\providecommand{\MRhref}[2]{%
  \href{http://www.ams.org/mathscinet-getitem?mr=#1}{#2}
}
\providecommand{\href}[2]{#2}

\end{document}